\documentclass[12pt,a4paper]{article}

\usepackage[top=3cm,bottom=3cm,left=2.4cm,right=2.4cm]{geometry}

\usepackage[utf8]{inputenc}
\usepackage[english]{babel}
\usepackage{amsmath,amssymb,amsfonts,amsthm,graphicx}
\usepackage{hyperref}

\begin{document}

\newcommand{\C}{\mathbb {C}}
\newcommand{\N}{\mathbb{N}}
\newcommand{\R}{\mathbb{R}}
\newcommand{\tf}{\mathcal{F}}

\renewcommand{\theenumi}{\roman{enumi}}
\renewcommand{\labelenumi}{\theenumi)}

\swapnumbers
\newtheorem{DEF-BFB}{Definition}[section]

\newtheorem{THM0.0}{Theorem}[section]
\newtheorem{DEF-STC}[THM0.0]{Definition}
\newtheorem{THM0.1}[THM0.0]{Theorem}
\newtheorem{THM0.2}[THM0.0]{Theorem}
\newtheorem{THM0.3}[THM0.0]{Theorem}

\newtheorem{3-EXPLE}{Example}[section]
\newtheorem{3-VDC}[3-EXPLE]{Theorem}
\newtheorem{3-REM1}[3-EXPLE]{Example}
\newtheorem{6-LEM}[3-EXPLE]{Lemma}
\newtheorem{3-REM2}[3-EXPLE]{Remark}
\newtheorem{3-VDC3}[3-EXPLE]{Theorem}
\newtheorem{3-VDC3bis}[3-EXPLE]{Remark}
\newtheorem{3-COR1}[3-EXPLE]{Theorem}
\newtheorem{3-REM2bis}[3-EXPLE]{Remark}
\newtheorem{3-OPTI1}[3-EXPLE]{Theorem}
\newtheorem{3-REM3}[3-EXPLE]{Remark}
\newtheorem{3-VDC4}[3-EXPLE]{Theorem}
\newtheorem{3-VDC4bis}[3-EXPLE]{Remark}
\newtheorem{3-OPTI2}[3-EXPLE]{Theorem}
\newtheorem{3-REM4}[3-EXPLE]{Remark}

\newtheorem{4-DEF1}{Definition}[section]
\newtheorem{4-REM8}[4-DEF1]{Remark}
\newtheorem{4-SCHRO1}[4-DEF1]{Theorem}
\newtheorem{4-SCHRO2}[4-DEF1]{Theorem}
\newtheorem{4-OPTI}[4-DEF1]{Corollary}
\newtheorem{4-REM9}[4-DEF1]{Remark}
\newtheorem{4-PREGEO1}[4-DEF1]{Theorem}
\newtheorem{4-PREGEO2}[4-DEF1]{Theorem}

\newtheorem{5-LEM3}{Lemma}[section]
\newtheorem{5-REM4}[5-LEM3]{Remark and Example}
\newtheorem{6-LEM2}[5-LEM3]{Lemma}
\newtheorem{5-SCHRO3}[5-LEM3]{Theorem}
\newtheorem{5-COR3}[5-LEM3]{Corollary}
\newtheorem{5-REM1}[5-LEM3]{Remark}

\title{Estimates of oscillatory integrals with stationary phase and singular amplitude:\\ Applications to propagation features for dispersive equations}

\author{Florent Dewez\footnote{Université de Valenciennes et du Hainaut-Cambrésis, LAMAV, FR CNRS 2956, Le Mont Houy, 59313 Valenciennes Cedex 9, France. E-mail: florent.dewez@univ-valenciennes.fr or florent.dewez@outlook.com}}

\date{}

\maketitle

\begin{abstract}
	In this paper, we study time-asymptotic propagation phenomena for a class of dispersive equations on the line by exploiting precise estimates of oscillatory integrals. We propose first an extension of the van der Corput Lemma to the case of phases which may have a stationary point of real order and amplitudes allowed to have an integrable singular point. The resulting estimates provide optimal decay rates which show explicitly the influence of these two particular points. Then we apply these abstract results to solution formulas of a class of dispersive equations on the line defined by Fourier multipliers. Under the hypothesis that the Fourier transform of the initial data has a compact support or an integrable singular point, we derive uniform estimates of the solutions in space-time cones, describing their motions when the time tends to infinity. The method permits also to show that symbols having a restricted growth at infinity may influence the dispersion of the solutions: we prove the existence of a cone, depending only on the symbol, in which the solution is time-asymptotically localized. This corresponds to an asymptotic version of the notion of causality for initial data without compact support.


\end{abstract}

\vspace{0.3cm}



\setcounter{section}{-1}
\section{Introduction}

\hspace{2.5ex} The time-asymptotic behaviour of solutions of dispersive equations depends strongly on the initial datum and certain phenomena may be exhibited by considering well-chosen sets of initial data. Following this idea, the aim of this paper is to exhibit propagation patterns for dispersive equations by choosing initial data in a \emph{bounded frequency band}, \emph{i.e.} their Fourier transform has a support contained in a bounded interval, or having a \textit{singular frequency}, \emph{i.e.} an integrable singular point of their Fourier transform. The method we propose is based on a careful application of precise estimates of oscillatory integrals with respect to a large parameter.

Such estimates can be obtained by employing the van der Corput Lemma \cite[Prop. 2, Chap. VIII]{stein} which exhibits in particular the decay rates of oscillatory integrals when the large parameter tends to infinity. Let us mention that this lemma and certain of its adaptations have been exploited in the literature to study solutions of evolution equations (see for example \cite{fam94}, \cite{benartzi1994} or \cite{msw}). In the present paper, we provide an extension of this lemma to the case of amplitude functions which may have an integrable singular point in the bounded integration interval in view of the above mentioned applications. Moreover, our extension covers also the case of phase functions having a stationary point of real order. The resulting estimates show in an explicit way the influence of the order of the stationary point and of the strength of the singular point on the decay.

These abstract results are then exploited to study the time-asymptotic behaviour of solutions of a family of dispersive equations defined by Fourier multipliers, permitting to cover especially equations of different nature like Schrödinger-type or hyperbolic equations. By applying our extension of the van der Corput Lemma to the Fourier solution formulas in well-chosen space-time cones, we show an asymptotic localization of the solutions produced by the restriction to bounded frequency bands, permitting especially to describe their time-asymptotic motions, and an effect of the singular frequency on the time-decay rate. We derive in particular an $L^{\infty}$-norm estimate with optimal decay rate for the solution of the free Schrödinger equation on the line with initial data having a singular frequency, completing the results obtained in \cite{fusion}.

Finally our method permits to show that, under certain restrictions on the growth of the symbol at infinity, the solution tends to be localized in a cone which depends on the symbol only. An illustration of this phenomenon is provided in the setting of the Klein-Gordon equation on the line: we prove that the solution tends to be localized in the light cone issued by the origin when the time tends to infinity, even in the case of initial data which are not compactly supported in space. This last result can be interpreted as an asymptotic version of the notion of causality for such initial data.

\section{Motivation and history} \label{hist-motiv}

\hspace{2.5ex} To explain the interest of considering initial data in bounded frequency bands or ha\-ving singular frequencies, let us give some physical explanations in the setting of the free Schrödinger equation on the line, namely,
\begin{equation} \label{eqschro}
	\left\{ \begin{array}{rl}
			& \hspace{-2mm} \big[ i \, \partial_t + \partial_{xx} \big] u_S(t) = 0 \\ [2mm]
			& \hspace{-2mm} u_S(0) = u_0
	\end{array} \right. \; ,
\end{equation}
where $t \geqslant 0$. In terms of Quantum Mechanics, the spatial position of a free particle with respect to time is described by the solution of the free Schrödinger equation \eqref{eqschro} under the condition $u_0 \in L^2(\R)$; in this case, the solution is called the \emph{wave function} of the particle. Moreover the solution is formally given by a wave packet of the form
\begin{equation} \label{solforschro}
	u_S(t,x) = \frac{1}{2 \pi} \int_{\R} \tf u_0 (p) \, e^{-itf_S(p) + ixp} \, dp \; ,
\end{equation}
where $f_S(p) := p^2$ and $\tf u_0$ is the Fourier transform of $u_0$. Now suppose that the initial datum $u_0$ is localized around a certain frequency $\overline{p} \in \R$, meaning that the particle at the time $t=0$ has a momentum localized around $\overline{p}$: for example, one may consider a function $u_0 \in \mathcal{S}(\R)$ such that its Fourier transform attains its maximum at $\overline{p}$ and has a support contained in $[\overline{p} - \varepsilon, \overline{p} + \varepsilon]$, where $\varepsilon > 0$ is small. By rewriting the expression of $f_S(p)$ as follows,
\begin{equation*}
	f_S(p) = f_S(\overline{p}) + f_S'(\overline{p})(p-\overline{p}) + \frac{1}{2} \, f_S''(\overline{p})(p-\overline{p})^2 = \overline{p}^2 + 2 \, \overline{p} (p-\overline{p}) + (p-\overline{p})^2 \; ,
\end{equation*}
one obtains 
\begin{equation*}
	u_S(t,x) = e^{-it \overline{p}^2 + ix\overline{p}} \, \frac{1}{2 \pi} \int_{\overline{p}-\varepsilon}^{\overline{p}+\varepsilon} \tf u_0 (p) \, e^{-i t \big( 2 \overline{p} (p-\overline{p}) + (p-\overline{p})^2\big) + ix(p-\overline{p})} \, dp =: e^{-it\overline{p}^2 + ix\overline{p}} \, E(t,x) \; ,
\end{equation*}
where $E(t,x)$ is actually called the \emph{envelope} of the wave packet: this factor contains the information on the position of the free particle. To describe roughly the motion of the envelope, one may proceed as follows: since the quadratic term $(p-\overline{p})^2$ is small as compared with $(p-\overline{p})$ for $p$ close to $\overline{p}$, it can be neglected in the expression of $E(t,x)$. Hence an approximation $\tilde{E}(t,x)$ of the envelope $E(t,x)$ is given by
\begin{equation*}
	\tilde{E}(t,x) = \frac{1}{2 \pi} \int_{\overline{p}-\varepsilon}^{\overline{p}+\varepsilon} \tf u_0 (p) \, e^{i (p-\overline{p}) ( -2 t \overline{p} + x)} \, dp \; .
\end{equation*}
Then we observe that this approximation is constant on the space-time half-line (or \emph{space-time direction}) $\frac{x}{t} = 2 \, \overline{p} = f_S'(\overline{p})$. Hence this approach indicates that the particle travels mainly at the speed $f_S'(\overline{p})$, which is actually called the \emph{group velocity} of the wave packet \eqref{solforschro}.\\

The above physical approach indicates that initial data having a certain localization in the frequency space may be useful to describe propagation features. Our aim is to establish rigorously such results by considering initial data which verify certain hypotheses on their Fourier transform:

\begin{DEF-BFB} \label{DEF-BFB}
	Let $p_1, p_2$ and $\tilde{p}$ be three finite real number such that $p_1 < p_2$.
	\begin{enumerate}
		\item A tempered distribution $u_0$ on $\R$ is in the \emph{bounded frequency band} $[p_1,p_2]$ if and only if its Fourier transform $\tf u_0$ is a complex-valued function such that
	\begin{equation*}
		supp \, \tf u_0 \subseteq [p_1,p_2] \; .
	\end{equation*}
		\item A tempered distribution $u_0$ on $\R$ has a \emph{singular frequency} at $\tilde{p}$ if and only if its Fourier transform $\tf u_0$ is a complex-valued function and $\tilde{p}$ is an integrable singular point of $\tf u_0$.
	\end{enumerate}	 
\end{DEF-BFB}

\noindent Let us explain such hypotheses and their consequences in the setting of the free Schrödinger equation. An initial datum $u_0 \in L^2(\R)$ in a bounded frequency band $[p_1,p_2]$ means that the particle at the time $t=0$ has a momentum localized in the interval $[p_1,p_2]$, but not necessarily around a particular value. According to the above approach, it does not seem possible to define in a rigorous way a group velocity for the resulting wave packet \eqref{solforschro} representing the solution of the free Schrödinger equation \eqref{eqschro}. However we observe that the wave packet has many frequency-components travelling at different speeds between $f'_S(p_1)$ and $f_S'(p_2)$, permitting to keep information on its motion: roughly speaking, one expects that the front of the wave packet moves at the speed of the fastest components while the back moves at the speed of the slowest ones, stretching out the envelope of the wave packet over time. In other words, one expects that the solution of equation \eqref{eqschro} with an initial datum in a frequency band $[p_1,p_2]$ is mainly localized in an interval of the form $\big[ f_S'(p_1) \, t, f_S'(p_2) \, t \big]$ for sufficiently large $t \geqslant 0$, describing the motion of the associated particle.

The presence of a singular frequency for an initial datum can be interpreted as a concentration of the momentum of the initial particle around this value. As above, defining in a precise way a group velocity for the resulting wave packet \eqref{solforschro} in the present case is not necessarily possible since the support of $\tf u_0$ is not required to be centered at the singular frequency or to be bounded. Nevertheless the strong accumulation of components travelling almost at the same speed is expected to affect the dispersion: these components need a longer time before being dissociated, diminishing in particular the time-decay rate of the solution. 

Let us also mention that the notion of bounded frequency band can be exploited to study the case of general initial data, that is to say initial data which are not necessarily in frequency bands. To do so, one can employ the following decomposition of a general initial datum $u_0$:
\begin{equation*}
	u_0 = \sum_{k \in K} \tf^{-1} \chi_{I_k} \, \tf u_0 =: \sum_{k \in K} u_{0,k} \; ,
\end{equation*}
where $K$ is a subset of $\mathbb{Z}$, $\{ I_k \}_{k \in K}$ is a family of bounded intervals such that
\begin{equation*}
	 \bigcup_{k \in K} I_k = \R\qquad , \qquad \forall \, k \neq l \quad I_k \cap I_l = \emptyset \; ,
\end{equation*}
and $\chi_{I_k}$ is the characteristic function of the interval $I_k$. Let us remark that the term $u_{0,k}$ is the component of $u_0$ in the bounded frequency band $I_k$. It turns out that the solution of the free Schrödinger equation \eqref{eqschro} with a general initial datum $u_0$ is actually the infinite sum of solutions of the same equation but with initial data given by $u_{0,k}$. Hence estimates of the solution can be derived from a precise study of each term of the infinite sum but, due to the superposition of all the bands $I_k$, the above mentioned localization of the solution for sufficiently large time may disappear.\\

In the paper \cite{fusion}, we established in a rigorous way results describing the motion of the solution of equation \eqref{eqschro} when the time tends to infinity. Our approach is based on the above physical ideas and so we considered initial data satisfying the following condition:\\

\noindent \textbf{Condition (C1$_{[p_1,p_2],\mu}$).} Fix $\mu \in (0,1]$ and let $p_1 < p_2$ be two finite real numbers.\\
A tempered distribution $u_0$ on $\R$ satisfies Condition \emph{(C1$_{[p_1,p_2],\mu}$)} if and only if $\tf u_0$ is a complex-valued function which verifies $supp \, \tf u_0 \subseteq [p_1,p_2]$ and
\begin{equation*} \label{singfreq}
	\forall \, p \in (p_1,p_2] \qquad \tf u_0 (p) = (p - p_1)^{\mu-1} \, \tilde{u}(p) \; ,
\end{equation*}
where $\tilde{u} \in \mathcal{C}^1\big( [p_1,p_2], \C \big)$ and $\tilde{u}(p_1) \neq 0$.\\

\noindent Let us note that an initial datum verifying this condition is in the frequency band $[p_1,p_2]$ and has a singular frequency at $p_1$. Inspired by the paper \cite{fam2012}, the method used in \cite{fusion} consists in expanding the solution formula \eqref{solforschro} to one term in certain \emph{space-time cones} of the form
\begin{equation} \label{defcone1}
	\mathfrak{C}(a,b) := \left\{ (t,x) \in (0,+\infty) \times \R \, \Big| \, 2 \, a < \frac{x}{t} < 2\, b \right\} \; ,
\end{equation}
by applying a stationary phase method. Hence we improved first the version of the stationary phase method of A. Erdélyi \cite[Section 9]{erdelyi} in order to obtain lossless error estimates for asymptotic expansions of oscillatory integrals. By applying then this refined method to the solution formula \eqref{solforschro}, we proved in Theorem 5.2 of \cite{fusion} that the preponderant time-decay rate of the solution is given by $t^{-\min\{\mu , \frac{1}{2}\}}$ inside cones of the type $\mathfrak{C}(p_1+ \varepsilon, p_2)$, where $\varepsilon > 0$ is sufficiently small, and by $t^{-\mu}$ in cones which are outside $\mathfrak{C}(p_1, p_2)$ according to Theorem 5.4. In particular in the $L^2$-case, namely when $\mu > \frac{1}{2}$, the decay outside the cone $\mathfrak{C}(p_1, p_2)$ is faster than inside, showing that the solution tends to be localized in the cone given by the frequency band when the time tends to infinity: this proves that the associated particle travels mainly in the interval $[2 p_1 \, t, 2 p_2 \, t] = \big[ f_S'(p_1) \, t, f_S'(p_2) \, t \big]$. Moreover Theorem 5.6 and Theorem 5.7 of \cite{fusion} emphasize the influence of the singular frequency on the decay: they show not only that the time-decay rate is diminished by the presence of the singular frequency $p_1$ but also that it seems to be slower in regions close to the space-time direction defined by $\frac{x}{t} = 2 \, p_1$. In Theorem 5.6, we expanded the solution to one term on $\frac{x}{t} = 2 \, p_1$, furnishing the slow time-decay rate $t^{-\frac{\mu}{2}}$; in Theorem 5.7, we established uniform and optimal estimates in curved space-time regions, depending on a parameter, along $\frac{x}{t} = 2 \, p_1$. The slow decay rate $t^{-\frac{\mu}{2}}$ is arbitrarily approached by a suitable choice of the parameter. The last results confirm the previous physical explanations: the time-decay rate is indeed diminished because of the concentration of components which travel at different speeds close to $2 \, p_1 = f_S'(p_1)$.

The method employed in \cite{fusion}, based on expansions to one term, is very precise locally in space but it does not permit to cover certain regions and hence to obtain uniform estimates in the whole space-time. This is due to the fact that the first terms and the remainders of the expansions given by the stationary phase method are not uniformly bounded with respect to the position of the stationary point in the case of singular amplitudes, whereas this is the case for the solution itself when $u_0$ satisfies Condition (C1$_{[p_1,p_2],\mu}$).

\section{Main results and references} \label{Sec-results}

\hspace{2.5ex} The aim of the present paper is to exhibit propagation features for certain dispersive equations, not by using asymptotic expansions as in \cite{fusion} but by creating a (rougher in a sense but closely related) method based on van der Corput type estimates, in other words estimates of the modulus of oscillatory integrals. This approach permits to avoid the blow-up of the expansions occurring in \cite{fusion} and hence to give uniform estimates in space-time cones as well as in the whole space-time, completing the results in \cite{fusion}. Moreover the explanations on wave packets given in the preceding section are not specific to the free Schrödinger equation and so, we propose in the present paper to consider the following class of dispersive equations on the line, including the Schrödinger setting:
\begin{equation} \label{eqdisp}
	\left\{ \begin{array}{rl}
			& \hspace{-2mm} \left[ i \, \partial_t - f \big(D \big) \right] u_f(t) = 0 \\ [2mm]
			& \hspace{-2mm} u_f(0) = u_0
	\end{array} \right. \; ,
\end{equation}
for $t \geqslant 0$, where the symbol $f$ of the Fourier multiplier $f(D)$ is supposed to satisfy $f'' > 0$. Note that equation \eqref{eqdisp} may not only describe physical phenomena, as in the Schrödinger setting, but also may be viewed as an intermediate equation which may appear when reducing higher order in time hyperbolic equations to first order equations (see for example Remark \ref{reduction}). \\


Before establishing results on the time-asymptotic behaviour of the solution of equation \eqref{eqdisp}, we study oscillatory integrals with respect to a large parameter $\omega$ of the type
\begin{equation} \label{oscillatory}
	\int_{p_1}^{p_2} U(p) \, e^{i \omega \psi(p)} \, dp \; ,
\end{equation}
in Section \ref{Sec-Theory}. The \emph{amplitude} $U : (p_1,p_2] \longrightarrow \C$ may be singular at $p_1$: we suppose that it is factorized as
\begin{equation} \label{amplitude}
	\forall \, p \in (p_1, p_2] \qquad U(p) = (p-p_1)^{\mu-1} \, \tilde{u}(p) \; ,
\end{equation}
where $\mu \in (0,1]$ and $\tilde{u} : [p_1,p_2] \longrightarrow \C$ is called the \emph{regular factor} of the amplitude. The \emph{phase function} $\psi : I \longrightarrow \R$, where $I$ is an open interval containing $[p_1,p_2]$, is allowed to have a stationary point $p_0$ of real order; more precisely, we suppose the factorization
\begin{equation} \label{phase}
	\forall \, p \in I \qquad \psi'(p) = | p - p_0 |^{\rho-1} \, \tilde{\psi}(p) \; ,
\end{equation}
where $\rho \in \R$ is larger than $1$ and $\tilde{\psi} : I \longrightarrow \R$, which satisfies $\big| \tilde{\psi} \big| > 0$ on $[p_1,p_2]$, is called the \emph{non-vanishing factor} of the phase. For example, smooth functions with vanishing first derivatives are included. The idea of supposing these factorizations, which are well suited for the formulation of the results of this paper, has been inspired by \cite{erdelyi}. \\
The first part of Section \ref{Sec-Theory} is devoted to the case of a phase function having a stationary point $p_0$ which is either inside or outside the interval of integration. In both cases, we furnish an estimate which is uniform with respect to the position of $p_0$, with explicit dependence of the decay rate on the order of the stationary point and on the strength of the singular point, thanks to a combination of the classical methods (see for example \cite{stein}) with the above well-adapted factorizations of the phase and of the amplitude. This provides an extension of the classical van der Corput Lemma:

\begin{THM0.0} \label{THM0.0}
	Let $\rho > 1$, $\mu \in (0,1]$ and choose $p_0 \in I$. Suppose that the functions $\psi : I \longrightarrow \R$ and $U : (p_1, p_2] \longrightarrow \C$ satisfy Assumption \emph{(P$_{p_0,\rho}$)} and Assumption \emph{(A$_{p_1,\mu}$)}, res\-pec\-ti\-ve\-ly. Moreover suppose that $\psi'$ is monotone on $\left\{ p \in I \, \big| \, p < p_0 \right\}$ and $\left\{ p \in I \, \big| \, p > p_0 \right\}$. Then we have
	\begin{equation*}
		\left| \int_{p_1}^{p_2} U(p) \, e^{i \omega \psi(p)} \, dp \, \right| \leqslant C(U,\psi) \, \omega^{-\frac{\mu}{\rho}} \; ,
	\end{equation*}
	for all $\omega > 0$, and the constant $C(U,\psi) > 0$ is given in the proof.
\end{THM0.0}

\noindent See Theorem \ref{3-COR1} for a complete statement. Assumption (A$_{p_1,\mu}$) and Assumption (P$_{p_0,\rho}$) are satisfied if and only if the functions $U$ and $\psi$ verify equalities \eqref{amplitude} and \eqref{phase} respectively, with additional hypotheses on the regularity of $\tilde{u}$ and $\tilde{\psi}$. \\
In the second part of Section \ref{Sec-Theory}, we suppose the absence of a stationary point inside $[p_1,p_2]$. Thus, if the phase has a stationary point, then it is outside the integration interval; in this case, the decay rate is better than the one obtained in Theorem \ref{THM0.0}, but the estimate is not uniform with respect to the position of the stationary point. In the applications, this result is essential to exhibit localization phenomena of the solution of \eqref{eqdisp} in space-time cones when the time tends to infinity.\\
Let us remark that all the decay rates provided in Section \ref{Sec-Theory} are proved to be optimal. This optimality is a consequence of \cite{fusion}, which gives asymptotic expansions to one term of the oscillatory integral.

In Section \ref{Sec-App1}, we study the time-asymptotic behaviour of solutions of equation \eqref{eqdisp}. To exploit the physical ideas given in the preceding section and deduce spatial information, we consider the Fourier solution formula:
\begin{equation} \label{solforeqdisp}
	u_f(t,x) = \frac{1}{2 \pi} \int_{\R} \tf u_0 (p) \, e^{-itf(p) + ixp} \, dp \; ,
\end{equation}
which permits the application of the results of Section \ref{Sec-Theory} thanks to a rewriting as an oscillatory integral with respect to time, inspired by \cite{fusion}.\\
To explain the results obtained in this section, let us give the definition of a \emph{space-time cone related to a symbol} $f$, extending the definition of the cone $\mathfrak{C}(a, b)$ given in \eqref{defcone1} and adapted to the study of the free Schrödinger equation:
\begin{DEF-STC} \label{DEF-STC}
	Let $a < b$ be two real numbers (eventually infinite) and let $f : \R \longrightarrow \R$ be a $\mathcal{C}^{\infty}$-function. We define the space-time cone $\mathfrak{C}_f(a, b)$ as follows:
	\begin{equation*}
		\mathfrak{C}_f(a, b) := \left\{ (t,x) \in (0,+\infty) \times \R \, \Big| \, f'(a) < \frac{x}{t} < f'(b) \right\} \; .
	\end{equation*}
	Let $\mathfrak{C}_f(a, b)^c$ be the complement of the cone $\mathfrak{C}_f(a, b)$ in $(0,+\infty) \times \R$ .
\end{DEF-STC}
\noindent We note in particular that $\mathfrak{C}(a, b) = \mathfrak{C}_{f_S}(a, b)$ where $f_S(p) = p^2$. Firstly we furnish uniform estimates of the solution \eqref{solforeqdisp} for initial data satisfying Condition (C1$_{[p_1,p_2],\mu}$) in arbitrary large space-time cones containing $\mathfrak{C}_f( p_1, p_2 )$ as well as in their complements:

\begin{THM0.1} \label{THM0.1}
	 Suppose that $u_0$ satisfies Condition \emph{(C1$_{[p_1,p_2],\mu}$)} and choose two finite real numbers $\tilde{p}_1 < \tilde{p}_2$ such that $\displaystyle [p_1,p_2] \subset (\tilde{p}_1, \tilde{p}_2) =: \tilde{I} $. Then we have
		\begin{align*}
			& \bullet \quad \forall \, (t,x) \in \mathfrak{C}_f(\tilde{p}_1,\tilde{p}_2) \qquad \big| u_f(t,x) \big| \leqslant c(u_0,f) \, t^{- \frac{\mu}{2}} \; ; \\[1mm]
			& \bullet \quad \forall \, (t,x) \in \mathfrak{C}_f(\tilde{p}_1,\tilde{p}_2)^c \qquad \big| u_f(t,x) \big| \leqslant c_{\tilde{I}}(u_0,f) \, t^{- \mu} \; .
		\end{align*}
		All the constants are given in the proof and the two decay rates are optimal.
\end{THM0.1}

\noindent See Theorem \ref{4-SCHRO1} for a complete statement. This result highlights the localization phenomenon produced by the bounded frequency band $[p_1,p_2]$ as explained in the preceding section in the Schrödinger case. In particular, it permits to derive an $L^{\infty}$-norm estimate of the solution for initial data satisfying Condition (C1$_{[p_1,p_2],\mu}$).\\
To study the influence of a singular frequency $p_1$ of order $\mu - 1$ in regions containing the space-time direction defined by $\frac{x}{t} = f'(p_1)$, we provide estimates of the solution in cones containing this direction as well as in cones which do not contain it:

\begin{THM0.2} \label{THM0.2}
	Suppose that $u_0$ satisfies Condition \emph{(C2$_{p_1, \mu}$)} (given in Section \ref{Sec-App1}), choose three finite real numbers $\varepsilon > 0$ and $\tilde{p}_1 < \tilde{p}_2$ such that $p_1 \notin [\tilde{p}_1, \tilde{p}_2]$. Then we have
		\begin{align*}
			& \bullet \quad \forall \, (t,x) \in \mathfrak{C}_f( p_1 - \varepsilon, p_1 + \varepsilon) \qquad \big| u_f(t,x) \big| \leqslant c^{(1)}(u_0,f) \, t^{-\frac{\mu}{2}} + c_{ \varepsilon}^{(2)}(u_0,f) \, t^{-1} \; ; \\[1mm]
			& \bullet \quad \forall \, (t,x) \in \mathfrak{C}_f(\tilde{p}_1, \tilde{p}_2 ) \qquad \big| u_f(t,x) \big| \leqslant c_{\tilde{p}_1, \tilde{p}_2}^{(1)}(u_0,f) \, t^{-\frac{1}{2}} + c_{\tilde{p}_1, \tilde{p}_2}^{(2)}(u_0,f) \, t^{-\mu} \\
			& \hspace{9cm} + c_{\tilde{p}_1, \tilde{p}_2}^{(3)}(u_0,f) \, t^{-1} \; .
		\end{align*}
		All the constants are given in the proof.
\end{THM0.2}

\noindent See Theorem \ref{4-PREGEO1} and Theorem \ref{4-PREGEO2} for a complete statement. This result proves that in cones containing the space-time direction given by the singular frequency, the influences on the decay rate of the singularity and the stationary point are combined. In all other cones, the decay rates given respectively by the singularity and the stationary point are in concurrence. Here Condition (C2$_{p_1, \mu}$) implies the presence of a singular frequency of order $\mu -1$ at $p_1$ but the initial datum is not necessarily in a  bounded frequency band; for example the support of $\tf u_0$ may be equal to an infinite interval. This permits to study the effect of the singular frequency without the influence of a bounded frequency band.

In Section \ref{Sec-App2}, we focus our attention on the effect of a symbol having a bounded first derivative on the dispersion of the solution. The interest of such symbols is not originated from a large family of physically interpretable equations of type \eqref{eqdisp} but rather from their usefulness to predict intrinsic time-asymptotic localization phenomena which may occur in the setting of higher order (in time) hyperbolic equations.\\
To explain this intrinsic localization, suppose that $f'(\R) = (a,b)$, where $a < b$ are two finite real numbers given by the limits of $f'$ at $-\infty$ and $+\infty$, and assume that the initial datum of equation \eqref{eqdisp} is localized around a frequency $\overline{p} \in \R$ as at the beginning of Section \ref{hist-motiv}. Then the resulting wave packet travels mainly at the group velocity $f'(\overline{p})$, which is always bounded from below by $a$ and from above by $b$. Consider now an initial datum in a frequency band $[p_1,p_2]$ whose associated solution is time-asymptotically localized in the space-time cone $\mathfrak{C}_f(p_1,p_2)$ according to Theorem \ref{THM0.1}. Under the same assumption $f'(\R) = (a,b)$, we have
\begin{equation*}
	\mathfrak{C}_f(p_1,p_2) \subset \left\{ (t,x) \in (0,+\infty) \times \R \, \Big| \, a < \frac{x}{t} < b \right\} =: \mathfrak{C}_f(-\infty,+\infty) \; ,
\end{equation*}
meaning that the solution of equation \eqref{eqdisp} can not be asymptotically localized outside the space-time cone $\mathfrak{C}_f(-\infty,+\infty)$, which depends on the symbol $f$ only, whatever the bounded frequency $[p_1,p_2]$ we choose.\\
Our aim in this section is to highlight this phenomenon by considering initial data which are not in a bounded frequency band. To do so, an explicit control of the behaviour of $f''$ at $\pm \infty$ is required for technical reasons:\\

\noindent \textbf{Condition (S$_{\beta_+, \beta_-,R}$).} Fix $\beta_- \geqslant \beta_+ > 1$ and $R \geqslant 1$.\\
A $\mathcal{C}^{\infty}$-function $f : \R \longrightarrow \R$ satisfies Condition (S$_{\beta_+, \beta_-,R}$) if and only if the second derivative of $f$ verifies
\begin{equation*}
	\exists \, c_+ \geqslant c_- > 0 \qquad \forall  \, |p| \geqslant R \qquad c_- \, |p|^{-\beta_-} \leqslant f''(p) \leqslant c_+ \, |p|^{-\beta_+} \; .
\end{equation*}
We note that the symbol $f_S$ related to the free Schrödinger equation \eqref{eqschro} does not satisfy the above condition, unlike the function $f_{KG}$ given by $f_{KG}(p) = \sqrt{c^4 + c^2 \, p^2 \,}$, where $c > 0$ is a constant. Under Condition (S$_{\beta_+, \beta_-,R}$), the first derivative of the symbol is bounded, leading to the existence of the space-time cone $\mathfrak{C}_f(-\infty,+\infty)$ in which the decay rate of the solution is slower than outside:

\begin{THM0.3} \label{THM0.3}
	Suppose that the symbol $f$ satisfies Condition \emph{(S$_{\beta_+, \beta_-,R}$)} and that $u_0$ satisfies Condition \emph{(C3$_{\mu,\alpha,r}$)} (given in Section \ref{Sec-App2}), where $\mu \in (0,1]$, $\alpha - \mu > \beta_-$ and $r \leqslant R$. Then we have
	\begin{align*}
		& \bullet \quad \forall \, (t,x) \in \mathfrak{C}_f(-\infty,+\infty) \qquad \big| u_f(t,x) \big| \leqslant c^{(1)}(u_0,f) \, t^{- \frac{\mu}{2}} + c^{(2)}(u_0,f) \, t^{- \frac{1}{2}} \; \\[1mm]
		& \bullet \quad \forall \, (t,x) \in \mathfrak{C}_f(-\infty,+\infty)^c \qquad \big| u_f(t,x) \big| \leqslant c_c^{(1)}(u_0,f) \, t^{- \mu} + c_c^{(2)}(u_0,f) \, t^{-1} \; .
	\end{align*}
	All the constants are given in the proof. The space-time cone $\mathfrak{C}_f(-\infty,+\infty)$ is defined by
	\begin{equation*}
		\mathfrak{C}_f(-\infty,+\infty) := \left\{ (t,x) \in (0,+\infty) \times \R \, \Big| \, \lim_{p \rightarrow - \infty} f'(p) < \frac{x}{t} < \lim_{p \rightarrow + \infty} f'(p) \right\} \; .
	\end{equation*}	
\end{THM0.3}

\noindent See Theorem \ref{5-SCHRO3} for a complete statement. Roughly speaking, an initial data satisfies Condition (C3$_{\mu,\alpha,r}$) if and only if it has a singular frequency of order $\mu - 1$ and $\tf u_0$ has a sufficient decay at infinity. To illustrate this intrinsic localization phenomenon in a simple setting, we consider the Klein-Gordon equation on the line. The solution of this hyperbolic equation of order 2 in time is actually the sum of the solutions of first order equations of type \eqref{eqdisp} with symbols $f_{KG}$ and $-f_{KG}$, where $f_{KG}$ is defined above. Hence Theorem \ref{THM0.3} is applicable to each term and we observe that the wave packets can not travel outside a light cone when the time tends to infinity. The appearance of this cone is closely related to the hyperbolic character of the Klein-Gordon equation.\\


Oscillatory integrals and their applications to evolution equations have been widely studied in the literature. One can mention \cite{reed-simon} in which the authors state the $L^1$-$L^{\infty}$ estimate for the unitary group generated by the free Hamiltonian for nonrelativistic quantum mechanics. Using interpolation, they obtain the well-known $L^p$-$L^q$ estimates for the free Schrödinger equation on $\R^n$. These results lead to Strichartz estimates which permit to study nonlinear variants of the equation.


The authors of \cite{msw} apply the van der Corput Lemma to the solution formulas of the wave equation and the Klein-Gordon equation on $\R^n$ to derive $L^{\infty}$-estimates.

A similar result was obtained in \cite{fam94} in the case of the Klein-Gordon equation on $\R$ with constant but different coefficients on the two half-axes. The author uses a spectral theoretic formula in order to apply the van der Corput Lemma to the solution of the equation.

In \cite{benartzi1994}, the authors consider a family of evolution equations given by \eqref{eqdisp}, where the symbols are of the form $f(p) = |p|^{\rho} + R(p)$, with $\rho \geqslant 2$ and $R : \R \longrightarrow \R$ is a regular function whose growth at infinity is controlled in a certain sense by $|p|^{\rho-1}$. They show that the operator $u_0 \in L^2(\R) \longmapsto u_f(t,.) \in L^2(\R)$ is unitary for all $t \geqslant 0$ and they establish the following estimate,
\begin{equation} \label{benartzi}
	\big\| u_f(t,.) \big\|_{L^{\infty}(\R)} \leqslant C \, t^{-\frac{1}{\rho}} \| u_0 \|_{L^1(\R)} \; ,
\end{equation}
for a certain constant $C > 0$, showing the dispersive nature of the equation. A Strichartz type estimate is then derived. The proof of estimate \eqref{benartzi} is based on the representation of the solution as a convolution of the initial datum with a distribution given by an oscillatory integral. By using a van der Corput type Lemma, estimates for this integral are obtained, leading to the result.\\
For comparison with the present paper, let us remark that our method based on Fourier solution formulas does not furnish a $L^1$-$L^{\infty}$ estimate but it permits to derive spatial information on the solution. This type of result does not seem possible with the method employed in \cite{benartzi1994}.

The time-decay rate of the free Schrö\-din\-ger equation is considered in \cite{cazenave1998} and \cite{cazenave2010}. In \cite{cazenave1998}, singular initial conditions are constructed to derive the exact $L^p$-time decay rates of the solution, which are slower than the classical results for regular initial conditions. In \cite{cazenave2010}, the authors construct initial conditions in Sobolev spaces (based on the Gaussian function), and they show that the related solutions has no definite $L^p$-time decay rates, nor coefficients, even though upper estimates for the decay rates are established. Both papers are based on special formulas for functions and their Fourier transform.\\
Though we do not furnish $L^p$-estimates for $p \neq \infty$ in the present paper, our method permits to cover a larger class of initial data, and hence to exhibit propagation patterns, and to consider more general symbols including $f_S(p) = p^2$.

In the setting of \cite{fam94}, the article \cite{fam2012} provides an asymptotic expansion to one term of the solution for initial data in bounded frequency bands, describing its time-asymptotic motion. Due to the potential step, there exist critical frequencies playing the same role as the singular frequencies introduced in the present paper. The paper \cite{fam2012}, which has inspired the study in \cite{fusion}, presents also a uniformity problem: the frequency band has to be chosen away from the critical frequencies, otherwise the expansion fails. The present paper solves this problem in the setting of \cite{fusion}, so we hope that our theoretical results may help to improve the understanding of the phenomena occuring in \cite{fam2012}.


One can also mention the results of \cite{fam0}, in which the authors consider the Schrödinger equation with sufficiently localized potential on a star-shaped network and provide $L^{\infty}$-decay estimates. A perturbation estimate shows that the solution is close to the free solution for initial data in high frequencies. In particular this result is applicable to the Schrödinger equation with potential on the line and permits to transfer some quantitative information from the free equation obtained in \cite{fusion} to the perturbed one.

To finish, we mention extensions of the van der Corput Lemma to the case of several integration variables which have been recently established to derive estimates of solutions of certain evolution equations on $\R^n$ for example.\\
In \cite{ruzhansky}, hypotheses on the radial behaviour of the phase in a neighbourhood of the stationary point permit to reduce the study to oscillatory integrals for one integration variable. By combining the standard calculations in the one-dimensional case and well-chosen assumptions on the phase and the amplitude, the author obtains the desired result. This approach permits to extend the notion of \emph{stationary point of integer order} to several variables and to provide estimates of oscillatory integrals with phases having this type of stationary point, but the resulting constants are not explicit.

Another example is given in \cite{zuily}, where a van der Corput-type estimate for several integration variables with explicit constant is established. To this end, the authors adapt the proof of the classical lemma to the case of several integration variables, leading to technical computations. Nevertheless their result is restricted to phases whose Hessian is supposed to be invertible, meaning that the order of the stationary point can not be larger than one in the one dimensional case.\\
Let us mention that the amplitudes in the two last papers are supposed to be smooth and compactly supported, meaning that they do not treat the case of singular amplitudes. An interesting outlook would be to find a suitable extension to several variables of the notion of \emph{singular point} for which van der Corput type estimates can be established.

\section{Stationary points of real order and singular amplitudes: van der Corput type estimates} \label{Sec-Theory}

\hspace{0.5cm} We start this section by stating the hypotheses on the phase function. Two examples are then given to illustrate these assumptions.\\

\noindent Let $p_1,p_2$ be two finite real numbers such that $p_1 < p_2$, and let $I$ be an open interval containing $[p_1,p_2]$.\\

\noindent \textbf{Assumption (P$_{p_0,\rho}$).} Let $p_0 \in I$ and $\rho > 1$.\\ A function $\psi : I \longrightarrow \R$ satifies Assumption (P$_{p_0,\rho}$) if and only if $\psi \in \mathcal{C}^1 \big( I \big) \cap \mathcal{C}^2 \big( I \backslash\{p_0\} \big)$ and there exists a function $\tilde{\psi} : I \longrightarrow \R$ such that
	\begin{equation*}
		\forall \, p \in I \qquad \psi'(p) = |p-p_0|^{\rho -1} \, \tilde{\psi}(p) \; ,
	\end{equation*}
	where $\big| \tilde{\psi} \big|: I \longrightarrow \R$ is assumed continuous and does not vanish on $I$.\\ The point $p_0$ is called \emph{stationary point} of $\psi$ of order $\rho -1$, and $\tilde{\psi}$ the \emph{non-vanishing factor} of $\psi$.\\
	
\noindent Let us comment on this choice. Firstly we consider the absolute value of $p-p_0$ because we want to include stationary points of non-integer order in the study. Secondly, the fact that $\tilde{\psi}$ does not vanish prevents this function from affecting the order of the stationary point $p_0$. Finally, the continuity of $\big| \tilde{\psi} \big|$ is sufficient to ensure the fact that $\displaystyle \min_{[p_1,p_2]} | \tilde{\psi} |$ exists and is non-zero; this quantity will be employed several times to establish the results of this section. Nevertheless we do not claim that we achieve full generality with these hypotheses. \\
Note that $\tilde{\psi}$ is actually continuously differentiable on $I \backslash \{p_0\}$, because
\begin{equation*}
	\forall \, p \neq p_0 \qquad \tilde{\psi}(p) = \frac{\psi'(p)}{|p-p_0|^{\rho-1}} \; .
\end{equation*}
This implies that $\tilde{\psi}$ has a constant sign on $\{ p \in I \, | \, p < p_0 \}$ and $\{ p \in I \, | \, p > p_0 \}$; note that the sign of $\tilde{\psi}$ can be different on each interval if $\tilde{\psi}$ has a discontinuity at the point $p_0$.\\
We illustrate the above Assumption (P$_{p_0,\rho}$) in the following two examples. In particular, the first example shows that smooth functions with vanishing first derivatives are included.

\begin{3-EXPLE} \label{3-EXPLE}
	\em \begin{enumerate}
		\item Let $\psi : I \longrightarrow \R$ be a function belonging to $\mathcal{C}^{N} \big( I \big)$ for a certain $N \geqslant 2$, and let $p_0 \in I$. Suppose that $\psi^{(k)}(p_0) = 0$ for $k=1, \ldots, N-1$. Then by Taylor's formula, we obtain
	\begin{align*}
		\psi'(p)	& = \frac{1}{(N-2)!} \int_{p_0}^p (p-x)^{N-2} \, \psi^{(N)}(x) \, dx \\
					& = \frac{(p-p_0)^{N-1}}{(N-2)!} \int_0^1 (1-y)^{N-2} \, \psi^{(N)} \big( y(p-p_0) + p_0 \big) \, dy \; ,
	\end{align*}
	for all $p \in I$. If we define $\tilde{\psi}$ as follows
	\begin{equation*}
		\tilde{\psi}(p) := \left\{ \begin{array}{rl}
			& \hspace{-4mm} \displaystyle \frac{1}{(N-2)!} \left(\frac{p-p_0}{|p-p_0|} \right)^{N-1} \int_0^1 (1-y)^{N-2} \, \psi^{(N)} \big( y(p-p_0) + p_0 \big) \, dy \; , \quad \text{if} \; p \neq p_0 \; , \\
			& \vspace{-0.3cm} \\
			& \hspace{-4mm} \displaystyle \frac{1}{(N-1)!} \, \psi^{(N)}(p_0) \; , \quad \text{if} \; p = p_0 \; ,
		\end{array} \right.
	\end{equation*}
	then $\displaystyle \psi'(p) = |p-p_0|^{N -1} \, \tilde{\psi}(p)$. Supposing $\big| \psi^{(N)} \big| > 0$ on $I$ implies that $\psi$ satisfies Assumption (P$_{p_0,N}$).
	\item Let $N \in \N$ such that $N \geqslant 2$ and choose $\alpha \in (N-1,N)$. Suppose that $\displaystyle \psi'(p) = |p|^{\alpha-1}$, for all $p \in \R$. In this case, $\psi \in \mathcal{C}^{N-1} \big( \R \big)$ but $\psi \notin \mathcal{C}^{N} \big( \R \big)$, and $\tilde{\psi} = 1$. Then Assumption (P$_{0,\alpha}$) is satisfied.
	\end{enumerate}
\end{3-EXPLE}

\vspace{4mm}

Now let us introduce the hypotheses concerning the amplitude function that we shall use throughout this section. \\

\noindent \textbf{Assumption (A$_{p_1,\mu}$).} Let $\mu \in (0,1]$.\\ A function $U : (p_1, p_2] \longrightarrow \C$ satisfies Assumption (A$_{p_1,\mu}$) if and only if there exists a function $\tilde{u} : [p_1,p_2] \longrightarrow \C$ such that
	\begin{equation*}
		\forall \, p \in (p_1,p_2] \qquad U(p )= (p - p_1)^{\mu -1} \, \tilde{u}(p) \; ,
	\end{equation*}
	where $\tilde{u}$ is assumed continuous on $[p_1,p_2]$, differentiable on $(p_1,p_2)$ with $\tilde{u}' \in L^1 (p_1,p_2)$, and $\tilde{u}(p_1) \neq 0$ if $\mu \neq 1$.\\ The point $p_1$ is called \emph{singular point} of $U$, and $\tilde{u}$ the \emph{regular factor} of $U$.\\
		
\noindent According to this assumption, the amplitude is singular at the left endpoint of the interval. The results of this section remain unchanged if we suppose that the singular point is at the right endpoint of the interval. Moreover in the case of an amplitude function which is singular inside the integration interval, the study can be reduced to the two preceding cases: it suffices to split the integral at the singular point.\\

Before providing the main results of this section, let us state a basic lemma which will be used several times.

\begin{6-LEM} \label{6-LEM}
	Let $\mu \in (0,1]$ and let $x,y \in \R_+$ such that $x \geqslant y$. Then we have
	\begin{equation*}
		x^{\mu} - y^{\mu} \leqslant (x-y)^{\mu} \; .
	\end{equation*}
\end{6-LEM}
\begin{proof}
	The case $\mu = 1$ is trivial so let us assume $\mu < 1$. If $y=0$ then the result is clear. Suppose $y \neq 0$, then the above inequality is equivalent to
	\begin{equation*}
		\left( \frac{x}{y} \right)^{\mu} - 1 \leqslant \left( \frac{x}{y} - 1 \right)^{\mu} \; .
	\end{equation*}
	Define the function $h : [1, +\infty) \longrightarrow \R$ by $h(t) := (t-1)^{\mu} - t^{\mu} + 1$. Then we note that for all $t > 1$,
	\begin{equation*}
		h'(t) = \mu \left( (t-1)^{\mu -1} - t^{\mu-1} \right) \geqslant 0 \; ,
	\end{equation*}
	since $\mu - 1 < 0$. It follows that $h(t) \geqslant h(1) = 0$, for all $t \in [1, +\infty)$, which proves the lemma.
\end{proof}

Our first aim in this section is to extend the classical van der Corput Lemma to the case of oscillatory integrals having a phase and an amplitude satisfying assumptions (P$_{p_0,\rho}$) and (A$_{p_1,\mu}$), respectively. This provides estimates of oscillatory integrals which are uniform with respect to the position of the stationary point $p_0$ and which exhibit the decay rates with respect to the large parameter.\\
In favour of readability, we prove two preliminary results before establishing our extension: in Theorem \ref{3-VDC}, the stationary point $p_0$ of order $\rho-1$ is supposed to belong to the integration interval $[p_1,p_2]$ while it is supposed to be outside this interval in Theorem \ref{3-VDC3}. In each theorem, the resulting estimate is uniform with respect to the position of $p_0$, the decay rate is given by $\omega^{-\frac{\mu}{\rho}}$ and an upper bound of the constant is given in terms of the regular factor $\tilde{u}$ of the amplitude and of the non-vanishing factor $\tilde{\psi}$ of the phase function. The combination of these two results leads to our extension of the van der Corput Lemma, stated in Theorem \ref{3-COR1}.

Let us now explain the main steps of the proof of Theorem \ref{3-VDC}; we shall follow similar steps to prove Theorem \ref{3-VDC3}. The proof is divided with respect to the size of the parameter $\omega$. In the case of small $\omega$, the integral can be estimated by the product of the $L^{\infty}$-norm of the regular part $\tilde{u}$ and the length of the interval to the power $\mu$, the exponent coming from the singular behaviour of the integrand at the point $p_1$ of order $\mu-1$. Then we exploit the fact that the interval is smaller than $\omega^{-\frac{1}{\rho}}$ to bound the integral by $\omega^{-\frac{\mu}{\rho}}$. For large $\omega$, we combine Stein's method \cite[Proposition 2, Chapter VIII]{stein}, which covers the case of stationary points of integer order and which is a generalization of Zygmund's method \cite[Lemma 4.3]{zygmund} for simple stationary points, and the above factorizations of $\psi'$ and $U$ inspired by \cite{erdelyi}. We decompose the integration interval in such way that $p_0$ and $p_1$ are contained in intervals whose length is proportional to $\omega^{-\frac{1}{\rho}}$. The integrals on these intervals are estimated by using the asymptotic smallness of their integration intervals. On the other intervals, we integrate by parts and we employ an upper bound for the amplitude as well as a lower bound for the first derivative of the phase, both bounds depending on $\omega$, to obtain the result. Let us note that we consider also the case of intermediate $\omega$; this situation can be studied by combining the methods used in the case of small and large $\omega$.

\begin{3-VDC} \label{3-VDC}
	Let $\rho > 1$, $\mu \in (0,1]$ and choose $p_0 \in [p_1,p_2]$. Suppose that the functions $\psi : I \longrightarrow \R$ and $U : (p_1, p_2] \longrightarrow \C$ satisfy Assumption \emph{(P$_{p_0,\rho}$)} and Assumption \emph{(A$_{p_1,\mu}$)}, respectively. Moreover suppose that $\psi'$ is monotone on $I_{p_0}^-$ and $I_{p_0}^+$, where
	\begin{equation*}
		I_{p_0}^- := \left\{ p \in I \, \big| \, p < p_0 \right\} \qquad , \qquad I_{p_0}^+ := \left\{ p \in I \, \big| \, p > p_0 \right\} \; .
	\end{equation*}
	Then we have
	\begin{equation*}
		\left| \int_{p_1}^{p_2} U(p) \, e^{i \omega \psi(p)} \, dp \right| \leqslant C(U,\psi) \, \omega^{-\frac{\mu}{\rho}} \; ,
	\end{equation*}
	for all $\omega > 0$, where the constant $C(U,\psi) > 0$ is given by
	\begin{equation*}
		C(U,\psi) := \frac{3}{\mu} \, \left\| \tilde{u} \right\|_{L^{\infty}(p_1,p_2)} + \left( 8 \left\| \tilde{u} \right\|_{L^{\infty}(p_1,p_2)} + 2 \left\| \tilde{u}' \right\|_{L^1(p_1,p_2)} \right) \left(\min_{p \in [p_1,p_2]} \left| \tilde{\psi}(p) \right| \right)^{-1} \; .
	\end{equation*}
\end{3-VDC}

\noindent Before proving this theorem, let us illustrate the monotonicity hypothesis on $\psi'$ by using the settings given in Example \ref{3-EXPLE}.

\begin{3-REM1}
	\em \begin{enumerate}
		\item In the setting of Example \ref{3-EXPLE} i), if $\big| \psi^{(N)} \big| > 0$ on $I$, then $\psi'$ is monotone on both intervals $I_{p_0}^-$ and $I_{p_0}^+$.\\
	Indeed if $N =2$, then it is clear that the hypothesis $\big| \psi'' \big| > 0$ implies the result. Suppose now that $N \geqslant 3$; then applying Taylor's formula to $\psi''$, namely
	\begin{equation*}
		\psi''(p) = \frac{1}{(N-3)!} \int_{p_0}^p (p-x)^{N-3} \, \psi^{(N)}(x) \, dx \; ,
	\end{equation*}
	for all $p \in I$, we observe that $\psi''$ has a constant sign on $I_{p_0}^-$ and $I_{p_0}^+$, which provides the result.
		\item In the setting of Example \ref{3-EXPLE} ii), we note that $\psi'$ is clearly monotone on $(-\infty,0)$ and on $(0,+\infty)$.
	\end{enumerate}	
\end{3-REM1}

\begin{proof}[Proof of Theorem \ref{3-VDC}]
	Let $p_0 \in (p_1,p_2)$ and let us suppose $\frac{p_0 - p_1}{2} \geqslant p_2 - p_0$ without loss of generality; the other case can be treated in a similar way. Note that we shall study the cases $p_0 = p_1$ and $p_0 = p_2$ at the end of the proof. Now let us divide the proof with respect to the size of $\omega$.
	\begin{itemize}
		\item \textit{Case $\omega > (p_2 - p_0)^{-\rho}$.} Define $\delta := \omega^{-\frac{1}{\rho}}$ and consider the following splitting of the integral:
		\begin{align*}
			\int_{p_1}^{p_2} U(p) \, e^{i \omega \psi(p)} \, dp	& \label{splitting} = \int_{p_1}^{p_1 + \delta} \dots \quad + \int_{p_1 + \delta}^{p_0-\delta} \dots \quad + \int_{p_0-\delta}^{p_0 + \delta} \dots \quad + \int_{p_0 + \delta}^{p_2} \dots \\
															& =: I^{(1)}(\omega) + I^{(2)}(\omega) + I^{(3)}(\omega) + I^{(4)}(\omega) \; . \nonumber
		\end{align*}
		Remark that this splitting is well-defined thanks to the hypothesis $\omega > (p_2 - p_0)^{-\rho}$. Let us estimate each integral.
		\begin{itemize}
			\item \textit{Study of $I^{(1)}(\omega)$.} We bound $I^{(1)}(\omega)$ in a simple way as follows:
			\begin{equation*}
				\Big| I^{(1)}(\omega) \Big| \leqslant \int_{p_1}^{p_1 + \delta} \left| U(p) \right| dp \leqslant \left\| \tilde{u} \right\|_{L^{\infty}(p_1,p_2)} \int_{p_1}^{p_1 + \delta} (p-p_1)^{\mu-1} \, dp = \frac{\left\| \tilde{u} \right\|_{L^{\infty}(p_1,p_2)}}{\mu} \, \delta^{\mu} \; .
			\end{equation*}
			\item \textit{Study of $I^{(2)}(\omega)$.} Here we shall suppose that $\tilde{\psi}$ is positive on $I_{p_0}^-$, which implies the positivity of $\psi'$; the case $\tilde{\psi} < 0$ can be studied in the same manner.  Since $\psi'$ does not vanish on $[p_1 + \delta, p_0 - \delta]$, the substitution $s = \psi(p)$ can be employed. Setting $\varphi := \psi^{-1}$, $s_1 := \psi(p_1 + \delta)$ and $s_2 := \psi(p_0 - \delta)$, we obtain
			\begin{align*}
				I^{(2)}(\omega)	& = \int_{s_1}^{s_2} U\big( \varphi(s) \big) \, \varphi'(s) \, e^{i \omega s} \, ds \\
							& = (i \omega)^{-1} \bigg( \Big[ (U \circ \varphi)(s) \,  \varphi'(s) \, e^{i \omega s} \Big]_{s_1}^{s_2} - \int_{s_1}^{s_2} \big( (U \circ \varphi) \, \varphi' \big)'(s) \, e^{i \omega s} \, ds \bigg) \; ;
			\end{align*}
			the last equality was obtained by integrating by parts.\\
			Now let us control the boundary terms and the integral. Firstly, we have
			\begin{equation} \label{estU}
				\big| U(p) \big| \leqslant \delta^{\mu-1} \left\| \tilde{u} \right\|_{L^{\infty}(p_1+\delta,p_0-\delta)} \leqslant \delta^{\mu-1} \left\| \tilde{u} \right\|_{L^{\infty}(p_1,p_2)} \; ,
			\end{equation}
			for all $p \in [p_1 + \delta, p_0 - \delta]$, since $U(p) = (p-p_1)^{\mu-1} \tilde{u}(p)$ by hypothesis. Moreover the fact that $\psi'$ satisfies Assumption (P$_{p_0,\rho}$) implies
			\begin{equation*}
				\forall \, p \in [p_1 + \delta, p_0 - \delta] \qquad \left| \psi'(p) \right| \geqslant \delta^{\rho-1} \, m \; ,
			\end{equation*}
			where $\displaystyle m := \min_{p \in [p_1,p_2]} \left| \tilde{\psi}(p) \right| > 0$. Combining this with the definition of $\varphi$ leads to
			\begin{equation} \label{estP}
				\forall \, s \in \left[s_1,s_2\right] \qquad \left| \varphi'(s) \right| \leqslant \delta^{1-\rho} \, m^{-1} \; .
			\end{equation}
			Inequalities \eqref{estU} and \eqref{estP} permit to estimate the boundary terms as follows,
			\begin{equation*}
				\bigg| \Big[ (U \circ \varphi)(s) \,  \varphi'(s) \, e^{i \omega s} \Big]_{s_1}^{s_2} \, \bigg| \leqslant 2 \, \left\| \tilde{u} \right\|_{L^{\infty}(p_1,p_2)} \, m^{-1} \, \delta^{\mu-\rho} \; .
			\end{equation*}
			It remains to control the integral. We have
			\begin{equation*}
				\big( (U \circ \varphi) \, \varphi' \big)' = (U' \circ \varphi) \, \left( \varphi' \right)^{\, 2} + (U \circ \varphi) \, \varphi'' \; ,
			\end{equation*}
			by the product rule; consequently,
			\begin{align}
				\bigg| \int_{s_1}^{s_2} \big( (U \circ \varphi) \, \varphi' \big)'(s) \, e^{i \omega s} \, ds \bigg|	& \leqslant \int_{s_1}^{s_2} \Big|(U' \circ \varphi)(s) \, \varphi'(s)^2 \Big| \, ds \nonumber \\
													& \quad \quad + \; \int_{s_1}^{s_2} \Big| (U \circ \varphi)(s) \, \varphi''(s) \Big| \, ds \nonumber \\
													& \leqslant \int_{s_1}^{s_2} \Big|(U' \circ \varphi)(s) \, \varphi'(s) \Big| \, ds \; \delta^{1-\rho} \, m^{-1} \nonumber \\
													& \quad \quad + \; \left\| U \right\|_{L^{\infty}(p_1+\delta,p_0-\delta)} \int_{s_1}^{s_2} \big| \varphi''(s) \big| \, ds \nonumber \\
													& \leqslant \int_{p_1+\delta}^{p_0-\delta} \big| U'(p) \big| \, dp \; \delta^{1-\rho} \, m^{-1} \nonumber \\
													& \label{estprincipal} \quad \quad + \; \delta^{\mu-1} \left\| \tilde{u} \right\|_{L^{\infty}(p_1,p_2)} \int_{s_1}^{s_2} \big| \varphi''(s) \big| \, ds \; .
			\end{align}
			The definition of $U$ implies
			\begin{align}
			\int_{p_1+\delta}^{p_0-\delta} \big| U'(p) \big| \, dp	& \leqslant \int_{p_1+\delta}^{p_0-\delta} \Big| (\mu-1) (p-p_1)^{\mu-2} \, \tilde{u}(p) \Big| \, dp \nonumber \\
	& \qquad \qquad + \; \int_{p_1+\delta}^{p_0-\delta} \Big| (p-p_1)^{\mu-1} \, \tilde{u}'(p) \Big| \, dp \nonumber \\
	& \leqslant \int_{p_1+\delta}^{p_0-\delta} (1-\mu) (p-p_1)^{\mu-2} \, dp \; \left\| \tilde{u} \right\|_{L^{\infty}(p_1,p_2)} \nonumber \\
	& \qquad \qquad + \; \delta^{\mu-1} \int_{p_1+\delta}^{p_0-\delta} \big| \tilde{u}'(p) \big| \, dp \nonumber \\
	& \label{est1} \leqslant \delta^{\mu-1} \left\| \tilde{u} \right\|_{L^{\infty}(p_1,p_2)} + \delta^{\mu-1} \left\| \tilde{u}' \right\|_{L^1(p_1,p_2)} \; ;
			\end{align}
			the last inequality was obtained employing the fact that
			\begin{equation*}
				\int_{p_1+\delta}^{p_0-\delta} (1-\mu) (p-p_1)^{\mu-2} \, dp = \delta^{\mu-1} - (p_0 - \delta - p_1)^{\mu-1} \leqslant \delta^{\mu-1} \; .
			\end{equation*}
			Moreover the relation $\displaystyle \varphi'' = \left( \frac{-\psi''}{\psi'^{\, 3}} \right) \circ \varphi$ provides the following equalities,
			\begin{equation*}
				\int_{s_1}^{s_2} \big| \varphi''(s) \big| \, ds = \int_{s_1}^{s_2} \left| \frac{-\psi''\big( \varphi(s)\big)}{\psi'\big( \varphi(s)\big)^3} \right| ds = \left| \int_{p_1 + \delta}^{p_0-\delta} \frac{-\psi''(p)}{\psi'(p)^2} \, dp \right| \; ,
			\end{equation*}
			the last equality comes from the change of variable $p = \varphi(s)$ and from the fact that $\psi''$ has a constant sign on $[p_1+\delta, p_0-\delta]$ thanks to the fact that $\psi'$ is monotonic on $I_{p_0}^-$. Then
			\begin{equation} \label{est2}
				\int_{s_1}^{s_2} \big| \varphi''(s) \big| \, ds = \left| \int_{p_1 + \delta}^{p_0-\delta} \left( \frac{1}{\psi'} \right)'(p) \, dp \right| = \left| \frac{1}{\psi'(p_0 - \delta)} - \frac{1}{\psi'(p_1 + \delta)} \right| \leqslant \delta^{1-\rho} m^{-1} \; ,
			\end{equation}
			where we used $| \psi'(p) | \geqslant \delta^{\rho-1} \, m$, for $p \in [p_1+\delta, p_0-\delta]$. Putting \eqref{est1} and \eqref{est2} in \eqref{estprincipal} provides
			\begin{equation*}
				\bigg| \int_{s_1}^{s_2} \big( (U \circ \varphi) \, \varphi' \big)'(s) \, e^{i \omega s} \, ds \bigg| \leqslant \left( 2 \left\| \tilde{u} \right\|_{L^{\infty}(p_1,p_2)} + \left\| \tilde{u}' \right\|_{L^1(p_1,p_2)} \right) m^{-1} \, \delta^{\mu-\rho} \; .
			\end{equation*}
			We are now able to estimate $I^{(2)}(\omega)$:
			\begin{equation*}
				\Big| I^{(2)}(\omega) \Big| \leqslant \left( 4 \left\| \tilde{u} \right\|_{L^{\infty}(p_1,p_2)} + \left\| \tilde{u}' \right\|_{L^1(p_1,p_2)} \right) m^{-1} \, \delta^{\mu-\rho} \, \omega^{-1} \; .
			\end{equation*}
			\item \textit{Study of $I^{(3)}(\omega)$.} As for $I^{(1)}(\omega)$, we bound the integral of $|U|$ on $[p_0 - \delta, p_0 + \delta]$ to provide an estimate of $I^{(3)}(\omega)$:
			\begin{equation*}
				\Big| I^{(3)}(\omega) \Big| \leqslant \frac{\left\| \tilde{u} \right\|_{L^{\infty}(p_1,p_2)}}{\mu} \, \big( (p_0 + \delta -p_1)^{\mu} - (p_0 - \delta - p_1)^{\mu} \big) \leqslant 2 \, \frac{\left\| \tilde{u} \right\|_{L^{\infty}(p_1,p_2)}}{\mu} \, \delta^{\mu} \; ,
			\end{equation*}
			where we applied Lemma \ref{6-LEM} to obtain the last inequality.
			\item \textit{Study of $I^{(4)}(\omega)$.} On $[p_0 + \delta, p_2]$, one can bound from below the absolute value of the first derivative of the phase function as follows,
			\begin{equation*}
				\big| \psi' \big| \geqslant \delta^{\rho-1} \, \min_{p \in [p_1,p_2]} \left| \tilde{\psi}(p) \right| = \delta^{\rho-1} \, m \; ,
			\end{equation*}
			and we have
			\begin{equation*}
				\forall \, p \in [p_0 + \delta, p_2] \qquad (p-p_1)^{\mu-1} \leqslant (p_0 + \delta - p_1)^{\mu-1} \leqslant \delta^{\mu-1} \; .
			\end{equation*}
			Following the lines of the study of $I^{(2)}(\omega)$ and using the two previous estimates, we obtain
			\begin{equation*}
				\Big| I^{(4)}(\omega) \Big| \leqslant \left( 4 \left\| \tilde{u} \right\|_{L^{\infty}(p_1,p_2)} + \left\| \tilde{u}' \right\|_{L^1(p_1,p_2)} \right) m^{-1} \, \delta^{\mu - \rho} \, \omega^{-1} \; .
			\end{equation*}
		\end{itemize}
		To conclude this first case, we replace $\delta$ by $\omega^{-\frac{1}{\rho}}$ leading to the desired estimate:
		\begin{align*}
			\left| \int_{p_1}^{p_2} U(p) \, e^{i \omega \psi(p)} \, dp \right|	& \leqslant \left| I^{(1)}(\omega) \right| + \left| I^{(2)}(\omega) \right| + \left| I^{(3)}(\omega) \right| + \left| I^{(4)}(\omega) \right| \\
		& \leqslant \frac{\left\| \tilde{u} \right\|_{L^{\infty}(p_1,p_2)}}{\mu} \, \omega^{-\frac{\mu}{\rho}} \; + \; 2 \, \frac{\left\| \tilde{u} \right\|_{L^{\infty}(p_1,p_2)}}{\mu} \, \omega^{-\frac{\mu}{\rho}} \\
		& \quad \qquad + \; 2 \left( 4 \left\| \tilde{u} \right\|_{L^{\infty}(p_1,p_2)} + \left\| \tilde{u}' \right\|_{L^1(p_1,p_2)} \right) m^{-1} \, \omega^{-\frac{\mu-\rho}{\rho}} \, \omega^{-1} \nonumber \\
		& =: C(U,\psi) \, \omega^{-\frac{\mu}{\rho}} \; ,
		\end{align*}
		where
		\begin{equation*}
			C(U,\psi) := \frac{3}{\mu} \, \left\| \tilde{u} \right\|_{L^{\infty}(p_1,p_2)} + \left( 8 \left\| \tilde{u} \right\|_{L^{\infty}(p_1,p_2)} + 2 \left\| \tilde{u}' \right\|_{L^1(p_1,p_2)} \right) m^{-1} \; .
		\end{equation*}
		\item \textit{Case $\big( \frac{p_0 - p_1}{2} \big)^{-\rho} < \omega \leqslant (p_2 - p_0)^{-\rho}$}. As above, we define $\delta := \omega^{-\frac{1}{\rho}}$ and we consider the following splitting of the integral:
			\begin{equation*}
				\int_{p_1}^{p_2} U(p) \, e^{i \omega \psi(p)} \, dp = \int_{p_1}^{p_1 + \delta} \dots \, + \, \int_{p_1 + \delta}^{p_0 - \delta} \dots \, + \, \int_{p_0 - \delta}^{p_0} \dots \, + \, \int_{p_0}^{p_2} \dots \; .
			\end{equation*}
			The three first integrals can be estimated using the methods of the first case, whereas the last integral can be controlled as follows,
			\begin{align}
				\left| \int_{p_0}^{p_2} U(p) \, e^{i \omega \psi(p)} \, dp \right|	& \leqslant \frac{\| \tilde{u} \|_{L^{\infty}(p_1,p_2)}}{\mu} \, \big( (p_2 - p_1)^{\mu} - (p_0 - p_1)^{\mu} \big) \nonumber \\
		&  \label{Int21} \leqslant \frac{\| \tilde{u} \|_{L^{\infty}(p_1,p_2)}}{\mu} \, (p_2 - p_0)^{\mu} \\
		& \label{Int2} \leqslant \frac{\| \tilde{u} \|_{L^{\infty}(p_1,p_2)}}{\mu} \, \delta^{\mu} \; ,
			\end{align}
			where we used Lemma \ref{6-LEM} to obtain inequality \eqref{Int21} and the fact that $p_2 - p_0 \leqslant \delta$ to establish inequality \eqref{Int2}. These arguments lead to
			\begin{align}
				\left| \int_{p_1}^{p_2} U(p) \, e^{i \omega \psi(p)} \, dp \right|	& \leqslant \left| \int_{p_1}^{p_1 + \delta} \dots \right| \, + \, \left| \int_{p_1 + \delta}^{p_0 - \delta} \dots \right| \, + \, \left| \int_{p_0 - \delta}^{p_0} \dots \right| \, + \, \left| \int_{p_0}^{p_2} \dots \right| \nonumber \\
				& \leqslant \frac{\| \tilde{u} \|_{L^{\infty}(p_1,p_2)}}{\mu} \, \delta^{\mu} + \left( 4 \left\| \tilde{u} \right\|_{L^{\infty}(p_1,p_2)} + \left\| \tilde{u}' \right\|_{L^1(p_1,p_2)} \right) m^{-1} \, \delta^{\mu-\rho} \, \omega^{-1} \nonumber \\
				& \label{largedelta2} \qquad + \frac{\| \tilde{u} \|_{L^{\infty}(p_1,p_2)}}{\mu} \, \delta^{\mu} + \frac{\| \tilde{u} \|_{L^{\infty}(p_1,p_2)}}{\mu} \, \delta^{\mu} \; ;
			\end{align}
			Replacing $\delta$ by $\omega^{-\frac{1}{\rho}}$ and observing that the constant which appears in \eqref{largedelta2} is smaller than $C(U, \psi)$ provides
			\begin{equation*}
				\left| \int_{p_1}^{p_2} U(p) \, e^{i \omega \psi(p)} \, dp \right| \leqslant C(U,\psi) \, \omega^{-\frac{\mu}{\rho}} \; .
			\end{equation*}
		\item \textit{Case $\omega \leqslant \big( \frac{p_0 - p_1}{2} \big)^{-\rho}$}. In this last case, we split the integral at the point $p_0$ and using the fact that $\omega \leqslant \big( \frac{p_0 - p_1}{2} \big)^{-\rho} \leqslant (p_2 - p_0)^{-\rho}$, we obtain
			\begin{align*}
				\left| \int_{p_1}^{p_2} U(p) \, e^{i \omega \psi(p)} \, dp \right|	& \leqslant \left| \int_{p_1}^{p_0} \dots \right| \, + \, \left| \int_{p_0}^{p_2} \dots \right| \\
						& \leqslant \frac{\| \tilde{u} \|_{L^{\infty}(p_1,p_2)}}{\mu} \, \big( (p_0 - p_1)^{\mu} + (p_2 - p_0)^{\mu} \big) \\
						& \leqslant 3 \, \frac{\| \tilde{u} \|_{L^{\infty}(p_1,p_2)}}{\mu} \, \omega^{-\frac{\mu}{\rho}} \; .
			\end{align*}
			We see that $C(U,\psi)$ is larger than the constant appearing in the right-hand side of the preceding inequality, leading to the result in this case.
	\end{itemize}
	Finally the desired estimate holds also for $p_0 = p_1$ and $p_0 = p_2$, since it is sufficient to adapt slightly the different splittings of the integral used in the present proof, and to carry out the same steps.
\end{proof}

\vspace{1mm}

\begin{3-REM2} \label{3-REM2} \em
	\begin{enumerate}
		\item The choice of the splitting points is optimized in view of the final decay rate. To prove that, we follow the indication given in the proof of Lemma 4.3 of \cite{zygmund}. Let us choose $\delta > 0$ sufficiently small to split the oscillatory integral as follows,
		\begin{equation*}
			\int_{p_1}^{p_2} U(p) \, e^{i \omega \psi(p)} \, dp = \int_{p_1}^{p_1 + \delta} \dots \quad + \int_{p_1 + \delta}^{p_0-\delta} \dots \quad + \int_{p_0-\delta}^{p_0 + \delta} \dots \quad + \int_{p_0 + \delta}^{p_2} \dots \; .
		\end{equation*}
		Applying the method employed in the case of large $\omega$ in the preceding proof gives an estimate of the form
		\begin{equation*}
			\left| \int_{p_1}^{p_2} U(p) \, e^{i \omega \psi(p)} \, dp \right| \leqslant f_{\omega}(\delta) \; ,
		\end{equation*}
		where $f_{\omega}(\delta) := c_1 \delta^{\mu} + c_2 \, \omega^{-1} \delta^{\mu-\rho}$, for certain constants $c_1 , c_2 > 0$. We note that $(f_{\omega})'$ vanishes at a unique point $\delta_0$ defined by
		\begin{equation*}
			\delta_0 := \left( \frac{\mu}{\rho-\mu} \, \frac{c_1}{c_2} \right)^{-\frac{1}{\rho}} \omega^{-\frac{1}{\rho}} \; .
		\end{equation*}
		Since $\displaystyle \lim_{\delta \rightarrow 0^+} f_{\omega}(\delta) = \lim_{\delta \rightarrow + \infty} f_{\omega}(\delta) = + \infty$, $\delta_0$ is then the minimum of $f_{\omega}$. Therefore the choice $\delta = \omega^{-\frac{1}{\rho}}$ is optimal regarding the decay rate.\\
		In particular, this splitting which depends on the parameter $\omega$ requires a decomposition of the proof with respect to the size of $\omega$. Indeed, the $\omega$-dependent cutting-points may leave the integration interval when $\omega$ is not sufficiently large.\\
		And we note that the constant $C(U,\psi)$ is surely not optimal since we do not choose exactly the minimum of $f_{\omega}$ for simplicity.
		\item Nevertheless the constant could be slightly improved in the case of regular amplitudes, namely $\mu = 1$ with $U = \tilde{u}$. Indeed, the study of $I^{(1)}(\omega)$ is not necessary in this situation and inequality \eqref{est1} can be simplified as follows,
		\begin{equation*}
			\int_{p_1}^{p_0-\delta} \big| U'(p) \big| \, dp \leqslant \left\| \tilde{u}' \right\|_{L^1(p_1,p_2)} \; .
		\end{equation*}
		It follows that we can estimate $I^{(2)}(\omega)$ and $I^{(4)}(\omega)$ more precisely,
		\begin{equation} \label{I2}
		\Big| I^{(j)}(\omega) \Big| \leqslant \left( 3 \left\| \tilde{u} \right\|_{L^{\infty}(p_1,p_2)} + \left\| \tilde{u}' \right\|_{L^1(p_1,p_2)} \right) m^{-1} \, \delta^{1 - \rho} \, \omega^{-1} \; ,
		\end{equation}
		with $j=2,4$, leading to
		\begin{equation*}
			C(U,\psi) := 2 \left\| \tilde{u} \right\|_{L^{\infty}(p_1,p_2)} + \left( 6 \left\| \tilde{u} \right\|_{L^{\infty}(p_1,p_2)} + 2 \left\| \tilde{u}' \right\|_{L^1(p_1,p_2)} \right) m^{-1} \; .
		\end{equation*}
		This refined constant will be used several times in Section \ref{Sec-App1}.
	\end{enumerate}
\end{3-REM2}

\vspace{4mm}

In Theorem \ref{3-VDC3}, we assume that the stationary point $p_0$ is outside the interval of integration $[p_1,p_2]$. In this case, the derivative of the phase function does not vanish inside the integration interval but it can be arbitrarily close to $0$ if the stationary point is close to this interval. The estimate that we provide does not depend on the position of the stationary point outside $[p_1,p_2]$, which makes that the resulting decay rate is the same as the one obtained in Theorem \ref{3-VDC}.\\
The proof of Theorem \ref{3-VDC3} is based on the same method as the one used to prove Theorem \ref{3-VDC}.

\begin{3-VDC3} \label{3-VDC3}
	Let $\rho > 1$, $\mu \in (0,1]$ and choose $p_0 \in I \backslash [p_1,p_2]$. Suppose that the functions $\psi : I \longrightarrow \R$ and $U : (p_1, p_2] \longrightarrow \C$ satisfy Assumption \emph{(P$_{p_0,\rho}$)} and Assumption \emph{(A$_{p_1,\mu}$)}, respectively. Moreover suppose that  $\psi'$ is monotone on $[p_1,p_2]$. Then we have
	\begin{equation*}
		\left| \int_{p_1}^{p_2} U(p) \, e^{i \omega \psi(p)} \, dp \right| \leqslant \tilde{C}(U,\psi) \, \omega^{-\frac{\mu}{\rho}} \; ,
	\end{equation*}
	for all $\omega > 0$, where the constant $\tilde{C}(U,\psi) > 0$ is given by
	\begin{equation*}
		\tilde{C}(U,\psi) := \frac{2}{\mu} \, \left\| \tilde{u} \right\|_{L^{\infty}(p_1,p_2)} + \left( 4 \left\| \tilde{u} \right\|_{L^{\infty}(p_1,p_2)} + \left\| \tilde{u}' \right\|_{L^1(p_1,p_2)} \right) \left(\min_{p \in [p_1,p_2]} \left| \tilde{\psi}(p) \right| \right)^{-1} \; .
	\end{equation*}
\end{3-VDC3}	

\begin{proof}
	We divide the proof with respect to the size of $\omega$.
	\begin{itemize}
	\item \textit{Case $\omega > \big( \frac{p_2 - p_1}{2}\big)^{-\rho}$}. We define $\delta := \omega^{-\frac{1}{\rho}}$ and we split the integral,
	\begin{align*}
		\int_{p_1}^{p_2} U(p) \, e^{i \omega \psi(p)} \, dp	& = \int_{p_1}^{p_1 + \delta} \dots \quad + \int_{p_1 + \delta}^{p_2-\delta} \dots \quad + \int_{p_2-\delta}^{p_2} \dots \\
															& =: \tilde{I}^{(1)}(\omega) + \tilde{I}^{(2)}(\omega) + \tilde{I}^{(3)}(\omega) \; ,
	\end{align*}
	where $\tilde{I}^{(1)}(\omega)$ and $\tilde{I}^{(3)}(\omega)$ are bounded from above by $\displaystyle \frac{\left\| \tilde{u} \right\|_{L^{\infty}(p_1,p_2)}}{\mu} \, \delta^{\mu}$. To estimate the integral $\tilde{I}^{(2)}(\omega)$, we follow the line of the method employed to study the integral $I^{(2)}(\omega)$ in the proof of Theorem \ref{3-VDC}, which provides
	\begin{equation*}
		\Big| \tilde{I}^{(2)}(\omega) \Big| \leqslant \left( 4 \left\| \tilde{u} \right\|_{L^{\infty}(p_1,p_2)} + \left\| \tilde{u}' \right\|_{L^1(p_1,p_2)} \right) m^{-1} \, \delta^{\mu - \rho} \, \omega^{-1} \; ;
	\end{equation*}
	we used the fact that
	\begin{equation*}
		\forall \, p \in [p_1 + \delta, p_2 - \delta] \qquad \big| U(p) \big| \leqslant \delta^{\mu-1} \| \tilde{u} \|_{L^{\infty}(p_1,p_2)} \; ,
	\end{equation*}
	and
	\begin{equation*}
		\forall \, p \in [p_1 + \delta, p_2 - \delta] \qquad \big| \psi'(p) \big| \geqslant \left\{ \begin{array}{rl}
			& \hspace{-3mm} (p_1 + \delta - p_0)^{\rho-1} \, m \geqslant \delta^{\rho-1} \, m \; , \quad \text{if } p_0 < p_1 \; , \\ [2mm]
			& \hspace{-3mm} (p_0 - p_2 + \delta)^{\rho-1} \, m \geqslant \delta^{\rho-1} \, m \; , \quad \text{if } p_0 > p_2 \; ,
	\end{array} \right.
	\end{equation*}
	with $\displaystyle m := \min_{p \in [p_1,p_2]} \left| \tilde{\psi} (p) \right|$. Finally we replace $\delta$ by $\omega^{-\frac{1}{\rho}}$ to conclude this case.
	\item \textit{Case $\omega \leqslant \big( \frac{p_2 - p_1}{2}\big)^{-\rho}$}. Here we have
	\begin{equation*}
		\left| \int_{p_1}^{p_2} U(p) \, e^{i \omega \psi(p)} \, dp \right| \leqslant \frac{\| \tilde{u} \|_{L^{\infty}(p_1,p_2)}}{\mu} (p_2 - p_1)^{\mu} \leqslant 2 \, \frac{\| \tilde{u} \|_{L^{\infty}(p_1,p_2)}}{\mu} \, \omega^{-\frac{\mu}{\rho}} \leqslant \tilde{C}(U,\psi) \, \omega^{-\frac{\mu}{\rho}} \; ,
	\end{equation*}
	which ends the proof.
	\end{itemize}
\end{proof}

\vspace{-1mm}

\begin{3-VDC3bis} \label{3-VDC3bis} \em
	In the case of regular amplitudes, one can use the estimate \eqref{I2} of $I^{(2)}(\omega)$ provided in Remark \ref{3-REM2}. In this situation, the constant $\tilde{C}(U,\psi)$ becomes
	\begin{equation*}
		\tilde{C}(U,\psi) := 2 \left\| \tilde{u} \right\|_{L^{\infty}(p_1,p_2)} + \left( 3 \left\| \tilde{u} \right\|_{L^{\infty}(p_1,p_2)} + \left\| \tilde{u}' \right\|_{L^1(p_1,p_2)} \right) \left(\min_{p \in [p_1,p_2]} \left| \tilde{\psi}(p) \right| \right)^{-1} \; .
	\end{equation*}
\end{3-VDC3bis}

\vspace{4mm}

Let us now state our extension of the van der Corput Lemma, which is actually a direct consequence of the two previous theorems.

\begin{3-COR1} \label{3-COR1}
	Let $\rho > 1$, $\mu \in (0,1]$ and choose $p_0 \in I$. Suppose that the functions $\psi : I \longrightarrow \R$ and $U : (p_1, p_2] \longrightarrow \C$ satisfy Assumption \emph{(P$_{p_0,\rho}$)} and Assumption \emph{(A$_{p_1,\mu}$)}, respectively. Moreover suppose that $\psi'$ is monotone on $I_{p_0}^-$ and $I_{p_0}^+$, where
	\begin{equation*}
		I_{p_0}^- := \left\{ p \in I \, \big| \, p < p_0 \right\} \qquad , \qquad I_{p_0}^+ := \left\{ p \in I \, \big| \, p > p_0 \right\} \; .
	\end{equation*}
	Then we have
	\begin{equation*}
		\left| \int_{p_1}^{p_2} U(p) \, e^{i \omega \psi(p)} \, dp \, \right| \leqslant C(U,\psi) \, \omega^{-\frac{\mu}{\rho}} \; ,
	\end{equation*}
	for all $\omega > 0$, where the constant $C(U,\psi) > 0$ is given by
	\begin{equation*}
		C(U,\psi) := \frac{3}{\mu} \, \left\| \tilde{u} \right\|_{L^{\infty}(p_1,p_2)} + \left( 8 \left\| \tilde{u} \right\|_{L^{\infty}(p_1,p_2)} + 2 \left\| \tilde{u}' \right\|_{L^1(p_1,p_2)} \right) \left(\min_{p \in [p_1,p_2]} \left| \tilde{\psi}(p) \right| \right)^{-1} \; .
	\end{equation*}
\end{3-COR1}

\begin{proof}
	Since the stationary point $p_0$ is allowed to be inside or outside the integration interval $[p_1,p_2]$, let us distinguish these two cases.
	\begin{itemize}
		\item \textit{Case $p_0 \in [p_1, p_2]$}. This corresponds to the setting of Theorem \ref{3-VDC} and so the integral is bounded by $C(U,\psi) \, \omega^{-\frac{\mu}{\rho}}$, where $C(U,\psi)$ is given in Theorem \ref{3-VDC}.
		\item \textit{Case $p_0 \notin [p_1, p_2]$}. In this case, either $[p_1,p_2] \subset I_{p_0}^-$ or $[p_1,p_2] \subset I_{p_0}^+$. Since $\psi'$ is assumed monotone on both intervals $I_{p_0}^-$ and $I_{p_0}^+$, Theorem \ref{3-VDC3} is applicable and then the integral is bounded by $\tilde{C}(U,\psi) \, \omega^{-\frac{\mu}{\rho}}$, where $\tilde{C}(U,\psi)$ is given in Theorem \ref{3-VDC3}.
	\end{itemize}
	Finally we remark that $\tilde{C}(U,\psi) \leqslant C(U,\psi)$, which concludes the proof.
\end{proof}

\vspace{1mm}

\begin{3-REM2bis} \label{3-REM2bis} \em
	As previously, we furnish a more precise constant in the case of regular amplitudes:
	\begin{equation*}
		C(U,\psi) := 2 \left\| \tilde{u} \right\|_{L^{\infty}(p_1,p_2)} + \left( 6 \left\| \tilde{u} \right\|_{L^{\infty}(p_1,p_2)} + 2 \left\| \tilde{u}' \right\|_{L^1(p_1,p_2)} \right) \left(\min_{p \in [p_1,p_2]} \left| \tilde{\psi}(p) \right| \right)^{-1} \; .
	\end{equation*}
\end{3-REM2bis}

\vspace{4mm}

In the following result, we prove the optimality of the decay rate given in Theorem \ref{3-COR1} under slightly stronger conditions. We show in fact that this decay rate is attained in the case of $p_0 = p_1$, where one can expect a superposition of the effects of the stationary point $p_0$ and the singular point $p_1$ of the amplitude.\\
Technically this result is based on an asymptotic expansion of the oscillatory integral to one term in this case. We use our recent results \cite[Theorems 2.3 and 2.7]{fusion} which are versions of the stationary phase method with explicit error estimates. Theorem 2.3 of \cite{fusion}, which covers the case of singular amplitudes, has been stated in \cite{erdelyi} with only rough indications of the steps of the proof. In \cite{fusion}, we have carried out all the details of this proof. Theorem 2.7 of \cite{fusion} is an improvement of the expansion result in \cite{erdelyi} for the case of regular amplitudes ($\mu = 1$).

\begin{3-OPTI1} \label{3-OPTI1}
	Suppose that the hypotheses of Theorem \ref{3-COR1} are satisfied. In addition to this, we assume that $\tilde{\psi}$ is right continuously differentiable at $p_0$ and $\tilde{u} \in \mathcal{C}^1\big( [p_1,p_2] \big)$ with $\tilde{u}(p_1) \neq 0$. \\
	Then the decay rate $\omega^{-\frac{\mu}{\rho}}$ given in Theorem \ref{3-COR1} is optimal and it is attained for $p_0 = p_1$.
\end{3-OPTI1}

\begin{proof}
	First of all, let us suppose that $p_0 = p_1$. Since the phase $\psi$ satisfies Assumption (P$_{p_1,\rho}$) in this case, the function $\tilde{\psi}$ has a constant sign on $(p_1,p_2]$ and it belongs to $\mathcal{C}^1\big( (p_1,p_2] \big)$. Hence the fact that $\tilde{\psi}$ is supposed to be right continuously differentiable at $p_0 = p_1$ implies that $\tilde{\psi}$ has a constant sign on $[p_1,p_2]$ and it belongs to $\mathcal{C}^1\big( [p_1,p_2] \big)$. \\
	Now let us suppose that $\tilde{\psi} > 0$ on $[p_1,p_2]$ without loss of generality. Hence the hypotheses of Theorem 2.3 of \cite{fusion} in the case $N=1$, $\rho_1 = \rho$, $\rho_2 = 1$, $\mu_1 = \mu$ and $\mu_2 = 1$ are satisfied and we obtain the following asymptotic expansion of the oscillatory integral with remainder estimates,
	\begin{equation*}
		\forall \, \omega > 0 \qquad \int_{p_1}^{p_2} U(p) \, e^{i \omega \psi(p)} \, dp = \sum_{j=1,2} \left( A_1^{(j)}(\omega) + R_1^{(j)}(\omega) \right) \; ,
	\end{equation*}
	where
	\begin{align*}
		& \displaystyle \bullet \quad A_1^{(1)}(\omega) := \frac{\rho^{\frac{\mu}{\rho}}}{\rho} \, \Gamma \Big( \frac{\mu}{\rho} \Big) \, e^{i \frac{\pi}{2} \, \frac{\mu}{\rho}} \, e^{i \omega \psi(p_1)} \, \frac{\tilde{u}(p_1)}{\tilde{\psi}(p_1)^{\frac{\mu}{\rho}}} \, \omega^{-\frac{\mu}{\rho}} \; , \\[2mm]
		& \displaystyle \bullet \quad A_1^{(2)}(\omega) := e^{-i \frac{\pi}{2}} \, e^{i \omega \psi(p_2)} \, \frac{U(p_2)}{\psi'(p_2)} \, \omega^{-1} \; , \\[2mm]
		& \displaystyle \bullet \quad \left| R_1^{(1)}(\omega) \right| \leqslant C^{(1)}(U,\psi, \nu) \, \omega^{-\frac{1}{\rho}} \; , \\[2mm]
		& \displaystyle \bullet \quad \left| R_1^{(2)}(\omega) \right| \leqslant C^{(2)}(U,\psi, \nu) \,  \omega^{-1} \; .
	\end{align*}
	The constants $C^{(1)}(U,\psi, \nu)$ and $C^{(2)}(U,\psi, \nu,)$ are independent from $\omega$ but both depend on a smooth cut-off function $\nu$ which separates the points $p_1$ and $p_2$. The above asymptotic expansion combined with the remainder estimates shows that $\omega^{-\frac{\mu}{\rho}}$ is the optimal decay rate.\\
	Let us remark that if $\mu = 1$, then Theorem 2.3 of \cite{fusion} gives the same decay rate for the first term $A_1^{(1)}(\omega)$ and for the remainder term $R_1^{(1)}(\omega)$, namely $\omega^{-\frac{1}{\rho}}$. To avoid this situation, one can employ Theorem 2.7 of \cite{fusion} which furnishes an estimate of $R_1^{(1)}(\omega)$ with a better decay rate than $\omega^{-\frac{1}{\rho}}$, and this fact assures that $\omega^{-\frac{1}{\rho}}$ is still the optimal decay rate for the oscillatory integral.
\end{proof}

\vspace{-1mm}

\begin{3-REM3}
	\emph{Theorem \ref{3-OPTI1} holds also when $\mu = 1$ and $\tilde{u}(p_1) = 0$, and in this case, the optimal decay rate $\omega^{-\frac{1}{\rho}}$ is attained for $p_0 = \tilde{p}$, if $\tilde{p} \in [p_1, p_2]$ satisfies $\tilde{u}(\tilde{p}) \neq 0$. To prove that, one can split the integral at $\tilde{p}$ and apply the stationary phase method to the two resulting integrals as in the preceding proof. We do not state this case in Theorem \ref{3-OPTI1} in favour of readability.}
\end{3-REM3}

\vspace{4mm}

Our second aim in the present section is to establish another estimate providing a faster decay rate of the oscillatory integral in the case of the absence of a stationary point inside the integration interval. More precisely we assume that the phase function is twice continuously differentiable and that its first derivative does not vanish on $[p_1,p_2]$. We obtain the decay rate $\omega^{-\mu}$ for singular amplitudes satisfying Assumption (A$_{p_1,\mu}$). This result will be necessary to exhibit localization phenomena for dispersive equations in Sections \ref{Sec-App1} and \ref{Sec-App2}.\\
When we want to apply the following result in the setting of Theorem \ref{3-VDC3} for comparison, we must suppose that the phase function is defined on an open interval $I$ which contains $[p_1,p_2]$ and that the stationary point $p_0$ of order $\rho-1$ of the phase belongs to $I \backslash [p_1,p_2]$, in other words it is outside the integration interval. In this case, Theorem \ref{3-VDC4} furnishes the better decay rate $\omega^{-\mu}$ as compared with the decay rate $\omega^{-\frac{\mu}{\rho}}$ given in Theorem \ref{3-VDC3}. Nevertheless the constant $C_c(U,\psi)$ of Theorem \ref{3-VDC4} tends to infinity when $p_0$ tends to $p_1$ or $p_2$, while the constant $\tilde{C}(U,\psi)$ provided in Theorem \ref{3-VDC3} is uniform with respect to the distance between $p_0$ and $[p_1,p_2]$.

\begin{3-VDC4} \label{3-VDC4}
	Let $\mu \in (0,1]$. Suppose that the function $U : (p_1,p_2] \longrightarrow \C$ satisfies Assumption \emph{(A$_{p_1,\mu}$)}. Moreover suppose that $\psi : I \longrightarrow \R$ belongs to $  \mathcal{C}^2\big([p_1,p_2]\big)$, and that $\psi'$ does not vanish and is monotone on $[p_1,p_2]$. Then we have
	\begin{equation*}
		\left| \int_{p_1}^{p_2} U(p) \, e^{i \omega \psi(p)} \, dp \right| \leqslant C_c(U,\psi) \, \omega^{-\mu} \; ,
	\end{equation*}
	for all $\omega > 0$, where the constant $C_c(U,\psi) > 0$ is given by
	\begin{equation*}
		C_c(U,\psi) := \frac{1}{\mu} \, \left\| \tilde{u} \right\|_{L^{\infty}(p_1,p_2)} + \left( 4 \left\| \tilde{u} \right\|_{L^{\infty}(p_1,p_2)} + \left\| \tilde{u}' \right\|_{L^1(p_1,p_2)} \right) \left( \min_{p \in [p_1,p_2]} \big| \psi'(p) \big| \right)^{-1} \; .
	\end{equation*}
\end{3-VDC4}

\begin{proof}
	We divide the proof with respect to $\omega$ one more time.
	\begin{itemize}
	\item \textit{Case $\omega > (p_2 - p_1)^{-1}$}. We define $\delta := \omega^{-1}$ and we consider the following splitting,
	\begin{align*}
		\int_{p_1}^{p_2} U(p) \, e^{i \omega \psi(p)} \, dp	& = \int_{p_1}^{p_1 + \delta} \dots \quad + \int_{p_1 + \delta}^{p_2} \dots \\
															& =: I_c^{(1)}(\omega) + I_c^{(2)}(\omega) \; .
	\end{align*}
	The integral $I_c^{(1)}(\omega)$ is bounded by $\displaystyle \frac{\left\| \tilde{u} \right\|_{L^{\infty}(p_1,p_2)}}{\mu} \, \delta^{\mu}$. Then we use the method employed to study the integral $I^{(2)}(\omega)$ in the proof Theorem \ref{3-VDC} in order to bound $I_c^{(2)}(\omega)$, since $\psi'$ does not vanish on $[p_1,p_2]$. But here, we bound $|\psi'|$ from below by $\displaystyle \min_{p \in [p_1,p_2]} | \psi'(p) | =: m > 0$, leading to
	\begin{equation*}
		\Big| I_c^{(2)}(\omega) \Big| \leqslant \left( 4 \left\| \tilde{u} \right\|_{L^{\infty}(p_1,p_2)} + \left\| \tilde{u}' \right\|_{L^1(p_1,p_2)} \right) m^{-1} \, \delta^{\mu - 1} \omega^{-1} \; .
	\end{equation*}
	Finally we replace $\delta$ by $\omega^{-1}$ to conclude this case.
	\item \textit{Case $\omega \leqslant (p_2 - p_1)^{-1}$}. We have
	\begin{equation*}
		\left| \int_{p_1}^{p_2} U(p) \, e^{i \omega \psi(p)} \, dp \right| \leqslant \frac{\| \tilde{u} \|_{L^{\infty}(p_1,p_2)}}{\mu} (p_2 - p_1)^{\mu} \leqslant \frac{\| \tilde{u} \|_{L^{\infty}(p_1,p_2)}}{\mu} \, \omega^{-\mu} \; ,
	\end{equation*}
	and we conclude the proof by noting that the constant which appears in the preceding inequality is smaller than $C_c(U,\psi)$.
	\end{itemize}
\end{proof}

\vspace{-1mm}

\begin{3-VDC4bis} \label{3-VDC4bis} \em
	Let us furnish a refinement of the constant $C_c(U,\psi)$ in the case of regular amplitudes. Here the integral $I_c^{(1)}(\omega)$ is not needed and according to Remark \ref{3-REM2}, the estimate of $I_c^{(2)}(\omega)$ is improvable. Then we obtain
		\begin{equation*}
			C_c(U,\psi) := \left( 3 \left\| \tilde{u} \right\|_{L^{\infty}(p_1,p_2)} + \left\| \tilde{u}' \right\|_{L^1(p_1,p_2)} \right) \left( \min_{p \in [p_1,p_2]} \big| \psi'(p) \big| \right)^{-1} \; .
		\end{equation*}
\end{3-VDC4bis}

\vspace{4mm}

In the last theorem of this section, we obtain the optimality of the decay rate given in Theorem \ref{3-VDC4} by applying Theorems 2.3 and 2.7 of \cite{fusion}, as we did in Theorem \ref{3-OPTI1}.

\begin{3-OPTI2} \label{3-OPTI2}
	Suppose that the hypotheses of Theorem \ref{3-VDC4} are satisfied. In addition to this, we assume that $\tilde{u} \in \mathcal{C}^1\big( [p_1,p_2] \big)$ with $\tilde{u}(p_1) \neq 0$. \\
	Then the decay rate $\omega^{-\mu}$ given in Theorem \ref{3-VDC4} is optimal.
\end{3-OPTI2}

\begin{proof}
	As in the proof of Theorem \ref{3-OPTI1}, we apply Theorem 2.3 of \cite{fusion} whose hypotheses are satisfied in the case $N=1$, $\rho_1 = \rho_2 = 1$, $\mu_1 = \mu$ and $\mu_2 = 1$. A new asymptotic expansion of the oscillatory integral with remainder estimates is then obtained,
	\begin{equation*}
		\forall \, \omega > 0 \qquad \int_{p_1}^{p_2} U(p) \, e^{i \omega \psi(p)} \, dp = \sum_{j=1,2} \left( \tilde{A}_1^{(j)}(\omega) + \tilde{R}_1^{(j)}(\omega) \right) \; ,
	\end{equation*}
	where
	\begin{align*}
		& \displaystyle \bullet \quad \tilde{A}_1^{(1)}(\omega) := \Gamma ( \mu ) \, e^{i \frac{\pi}{2} \, \mu} \, e^{i \omega \psi(p_1)} \, \frac{\tilde{u}(p_1)}{\psi'(p_1)^{\mu}} \, \omega^{-\mu} \; , \\[2mm]
		& \displaystyle \bullet \quad \tilde{A}_1^{(2)}(\omega) := e^{-i \frac{\pi}{2}} \, e^{i \omega \psi(p_2)} \, \frac{U(p_2)}{\psi'(p_2)} \, \omega^{-1} \; , \\[2mm]
		& \displaystyle \bullet \quad \left| \tilde{R}_1^{(1)}(\omega) \right| \leqslant \tilde{C}^{(1)}(U,\psi, \nu) \, \omega^{-1} \; , \\[2mm]
		& \displaystyle \bullet \quad \left| \tilde{R}_1^{(2)}(\omega) \right| \leqslant \tilde{C}^{(2)}(U,\psi, \nu) \,  \omega^{-1} \; .
	\end{align*}
	As above, $\nu$ is a smooth cut-off function separating the points $p_1$ and $p_2$, and the constants $\tilde{C}^{(1)}(U,\psi, \nu)$ and $\tilde{C}^{(2)}(U,\psi, \nu,)$ are independent from $\omega$. Hence we can conclude that $\omega^{-\mu}$ is the optimal decay rate.\\
	And if $\mu = 1$, then we employ Theorem 2.7 of \cite{fusion} to obtain more precise estimates for $\tilde{R}_1^{(1)}(\omega)$ and $\tilde{R}_1^{(2)}(\omega)$, furnishing better decay rates than $\omega^{-1}$ for these remainder terms, and so $\omega^{-1}$ is still the optimal decay rate.
\end{proof}

\vspace{-1mm}

\begin{3-REM4}
	\emph{When $\mu = 1$ and $\tilde{u}(p_1) = \tilde{u}(p_2) = 0$, the decay rate may be faster than $\omega^{-1}$. In this case, it depends on the regularity of $\psi$ and $U$, and on the values of the successive derivatives of these functions at the endpoints of the integration interval (see \cite[p. 331]{stein}).}
\end{3-REM4}

\section{Applications to a class of dispersive equations: influence of bounded frequency bands and singular frequencies on dispersion} \label{Sec-App1}

\hspace{0.5cm} In this section, we consider the class of evolution equations on the line defined by Fourier multipliers whose symbol has a positive second derivative, and we suppose that the initial data are in a bounded frequency band or have a singular frequency. Our aim is to exhibit propagation patterns produced by the frequency band and by the singular frequency by estimating the solutions in space-time cones thanks to the preceding abstract estimates of oscillatory integrals. Moreover studying this class of first order (in time) equations may be useful when considering higher order hyperbolic equations (see Theorem \ref{5-COR3} or Remark \ref{reduction}). \\

Let us now describe the setting of the present section: let $f : \R \longrightarrow \R$ be a $\mathcal{C}^{\infty}$-function such that all derivatives grow at most as a polynomial at infinity. We can associate with such a \textit{symbol} $f$ an operator $f(D) : \mathcal{S}(\R) \longrightarrow \mathcal{S}(\R)$ defined by
\begin{equation*}
	\forall \, x \in \R \qquad f(D) u(x) := \frac{1}{2 \pi} \int_{\R} f(p) \, \tf u(p) \, e^{i x p} \, dp = \tf^{-1} \Big( f \, \tf u \Big)(x) \; ,
\end{equation*}
where $\tf u$ is the Fourier transform of $u \in \mathcal{S}(\R)$, namely $\displaystyle \tf u(p) = \int_{\R} u(x) \, e^{- i x p} \, dx$. Since all the derivatives of the symbol $f$ grow at most as a polynomial at infinity, $f(D)$ can be extended to a map from the tempered distributions $\mathcal{S}'(\R)$ to itself.  The operator $f(D) : \mathcal{S}'(\R) \longrightarrow \mathcal{S}'(\R)$ is called a \textit{Fourier multiplier}.\\
Then for such an operator, we can introduce the following evolution equation on the line,
\begin{equation} \label{evoleq}
	\left\{ \begin{array}{rl}
			& \hspace{-2mm} \left[ i \, \partial_t - f \big(D\big) \right] u_f(t) = 0 \\ [2mm]
			& \hspace{-2mm} u_f(0) = u_0
	\end{array} \right. \; ,
\end{equation}
for $t \geqslant 0$. Throughout this section, we shall suppose that $f'' > 0$; an important example of such an equation is given by the free Schrödinger equation whose symbol is $f_S(p) = p^2$. Let us remark that one can also establish similar results to those of the present section when the second derivative of the symbol is supposed to be negative. Supposing $u_0 \in \mathcal{S}'(\R)$, the equation \eqref{evoleq} has a unique solution in $\displaystyle \mathcal{C}^1\big( \R_+ , \mathcal{S}'(\R) \big)$, given by the following solution formula,
\begin{equation} \label{formula}
	u_f(t) = \tf^{-1} \Big( e^{-i t f} \tf u_0 \Big) \; .
\end{equation}

Now we recall the definition of a space-time cone related to a symbol $f$.

\begin{4-DEF1}
	Let $a < b$ be two real numbers (eventually infinite) and let $f : \R \longrightarrow \R$ be a $\mathcal{C}^{\infty}$-function. We define the space-time cone $\mathfrak{C}_f(a, b)$ as follows:
	\begin{equation*}
		\mathfrak{C}_f(a, b) := \left\{ (t,x) \in (0,+\infty) \times \R \, \Big| \, f'(a) < \frac{x}{t} < f'(b) \right\} \; .
	\end{equation*}
	Let $\mathfrak{C}_f(a, b)^c$ be the complement of the cone $\mathfrak{C}_f(a, b)$ in $(0,+\infty) \times \R$ .
\end{4-DEF1}

\vspace{4mm}

In the first result of this section, we show that the solution of the equation \eqref{evoleq} for initial data in a bounded frequency band $[p_1,p_2]$ tends to be localized in the space-time cone $\mathfrak{C}_f(p_1, p_2)$ when the time tends to infinity. To do so, we furnish estimates with optimal decay rates of the solution inside arbitrary space-time cones containing $\mathfrak{C}_f(p_1, p_2)$ as well as in their complements. It turns out that the resulting decay rates are always slower inside the cones than outside, proving the above mentioned time-asymptotic localization. Moreover we allow the point $p_1$ to be a singular frequency: the resulting time-decay rates are then slower than they are in the regular case. In \cite[Theorems 5.2 and 5.4]{fusion}, different decay rates have already been obtained in the Schrödinger case by expanding the solution in certain space-time cones, but without uniformity as explained in the introduction of the present paper.\\
The first step of the proof consists in rewriting the solution formula as an oscillatory integral of the form \eqref{oscillatory}. In particular, the resulting phase function depends explicitly on the parameters $x$ and $t$ and has at most one stationary point which depends on the quotient $\frac{x}{t}$. The following step is to apply the results of the preceding section. To do so, we divide the proof with respect to the value of $\frac{x}{t}$: if $\frac{x}{t}$ is in a neighborhood of the integration interval, then we apply Theorem \ref{3-COR1} leading to a uniform estimate in a cone containing $\mathfrak{C}_f(p_1, p_2)$ with the slow decay $t^{-\frac{\mu}{2}}$. Otherwise, we obtain the better decay rate $t^{-\mu}$ outside the cone by applying Theorem \ref{3-VDC4}. The optimality of the rates is a direct consequence of Theorem \ref{3-OPTI1} and Theorem \ref{3-OPTI2}. \\

\noindent \textbf{Condition (C1$_{[p_1,p_2],\mu}$).} Fix $\mu \in (0,1]$ and let $p_1 < p_2$ be two finite real numbers. \\
	A tempered distribution $u_0$ on $\R$ satisfies Condition (C1$_{[p_1,p_2],\mu}$) if and only if $supp \, \tf u_0 \subseteq [p_1,p_2]$ and $\tf u_0$ verifies Assumption (A$_{p_1,\mu}$) on $[p_1,p_2]$, where the regular factor $\tilde{u}$ is supposed to belong to $\mathcal{C}^1\big( [p_1,p_2] \big)$ and $\tilde{u}(p_1) \neq 0$.

\vspace{2mm}

\begin{4-REM8} \label{4-REM8} \em
	\begin{enumerate}
		\item The subset of tempered distributions satisfying Condition (C1$_{[p_1,p_2],\mu}$) is non-empty. Indeed if a function $U$ verifies Assumption (A$_{p_1,\mu}$) with $supp \, U \subseteq [p_1,p_2]$ and with a regular factor belonging to $\mathcal{C}^1 \big([p_1,p_2] \big)$, then $U$ is an integrable function and so it belongs to $\mathcal{S}'(\R)$. Since the Fourier transform is a bijection on $\mathcal{S}'(\R)$, there exists $u_0 \in \mathcal{S}'(\R)$ such that $U = \tf u_0$, and hence $u_0$ satisfies Condition (C1$_{[p_1,p_2],\mu}$).
		\item Since the support of $\tf u_0$ is contained in a bounded interval, $u_0$ is in fact an analytic function on $\R$.
		\item Condition (C1$_{[p_1,p_2],\mu}$) implies that the initial condition has a singular frequency at the left endpoint of its bounded frequency band. As explained just after the statement of Assumption (A$_{p_1,\mu}$) in Section \ref{Sec-Theory}, the result in the case of a singular frequency inside the frequency band is analogous to the result stated in Theorem \ref{4-SCHRO1}.
		\item Thanks to the integrability of $\tf u_0$, the solution formula given in \eqref{formula} defines a complex-valued function on $\R_+ \times \R$ as follows,
		\begin{equation} \label{formula2}
			\forall \, (t,x) \in \R_+ \times \R \qquad u_f(t,x) = \frac{1}{2 \pi} \int_{\R} \tf u_0(p) \, e^{-i t f(p) + i x p} \, dp \; .
		\end{equation}
	\end{enumerate}
\end{4-REM8}

\vspace{1mm}

\begin{4-SCHRO1} \label{4-SCHRO1}
	 Suppose that $u_0$ satisfies Condition \emph{(C1$_{[p_1,p_2],\mu}$)} and choose two finite real numbers $\tilde{p}_1 < \tilde{p}_2$ such that $\displaystyle [p_1,p_2] \subset (\tilde{p}_1, \tilde{p}_2) =: \tilde{I} $. Then we have
		\begin{equation*}
			\forall \, (t,x) \in \mathfrak{C}_f(\tilde{p}_1,\tilde{p}_2) \qquad \big| u_f(t,x) \big| \leqslant c(u_0,f) \, t^{- \frac{\mu}{2}} \; ,
		\end{equation*}
		where the constant $c(u_0,f) > 0$ is given by \eqref{c_I}. Moreover we have
		\begin{equation*}
			\forall \, (t,x) \in \mathfrak{C}_f(\tilde{p}_1,\tilde{p}_2)^c \qquad \big| u_f(t,x) \big| \leqslant c_{\tilde{I}}(u_0,f) \, t^{- \mu} \; ,
		\end{equation*}
		where the constant $c_{\tilde{I}}(u_0,f) > 0$ is given by \eqref{c_I^c}. And the two decay rates are optimal.
\end{4-SCHRO1}

\begin{proof}
	We consider the solution formula given by \eqref{formula2} and we factorize the phase function $p \longmapsto xp - tf(p)$ by $t$, which gives
	\begin{equation*}
		\forall \, (t,x) \in (0,+\infty) \times \R \qquad u_f(t,x) = \int_{p_1}^{p_2} U(p) \, e^{i t \psi(p)} \, dp \; ,
	\end{equation*}
	where
	\begin{equation*}
		\left\{ \begin{array}{rl}
				& \displaystyle \forall \, p \in (p_1,p_2] \qquad U(p) := \frac{1}{2\pi}  \, \tf u_0(p) = \frac{1}{2\pi}  \, (p-p_1)^{\mu-1} \, \tilde{u}(p) \; , \\
				& \vspace{-0.3cm} \\
				& \displaystyle \forall \, p \in \R \qquad \psi(p) := \frac{x}{t} \, p - f(p) \; .
		\end{array} \right.
	\end{equation*}
	By hypothesis, the function $U$ verifies Assumption (A$_{p_1,\mu}$) on $[p_1,p_2]$. Moreover, we recall that $f''$ is supposed to be positive on $\R$, which implies that $f' : \R \longrightarrow f'(\R)$ is strictly increasing. It follows that the function $\psi'$ given by
	\begin{equation*}
		\forall \, p \in \R \qquad \psi'(p) = \frac{x}{t} - f'(p) \; ,
	\end{equation*}
	is strictly decreasing on $\R$. In particular, if a stationary point $p_0$ exists then it is unique and it is defined by
	\begin{equation*}
		p_0 = \big(f'\big)^{-1} \Big(\frac{x}{t} \Big) \; .
	\end{equation*}
	Hence the existence of a stationary point as well as its position with respect to the integration interval depends on the value of $\frac{x}{t}$. This leads us to divide the rest of the proof into two parts.
	\begin{enumerate}
		\item \textit{Case $\frac{x}{t} \in f'\big(\tilde{I}\big)$.} In this case, the stationary point $p_0$ exists and it belongs to $\tilde{I} := (\tilde{p}_1, \tilde{p}_2)$. Moreover the fact that $\psi'' = - f'' < 0$ implies $\psi''(p_0) \neq 0$. Consequently, according to Example \ref{3-EXPLE} i), the function $\psi : \R \longrightarrow \R$ satisfies Assumption (P$_{p_0,2}$) with
		\begin{equation} \label{regular-phase}
			\tilde{\psi}(p) = \left\{ \begin{array}{rl}
				& \displaystyle \frac{p-p_0}{|p-p_0|} \int_0^1 -f''\big(y(p-p_0) + p_0\big) \, dy \; , \quad \text{if} \; p \neq p_0 \; , \\
				& \vspace{-0.3cm} \\
				& \displaystyle -f''(p_0) \; , \quad \text{if} \; p = p_0 \; ,
			\end{array} \right.
		\end{equation}
		and $\big| \tilde{\psi}(p) \big| \geqslant m > 0$ for all $p \in [p_1,p_2]$, where $\displaystyle m := \min_{p \in [p_1,p_2]} f''(p) > 0$. So we can apply Theorem \ref{3-COR1} with $\rho = 2$, which gives
		\begin{equation*}
			\forall \, (t,x) \in \mathfrak{C}_f(\tilde{p}_1,\tilde{p}_2) \qquad \big| u_f(t,x) \big| = \left| \int_{p_1}^{p_2} U(p) e^{it \psi(p)} \, dp \right| \leqslant c(u_0,f) \, t^{-\frac{\mu}{2}} \; ,
		\end{equation*}
		where
		\begin{equation} \label{c_I}
			c(u_0,f) := \frac{1}{2\pi} \, \frac{3}{\mu} \, \left\| \tilde{u} \right\|_{L^{\infty}(p_1,p_2)} + \frac{1}{\pi} \left( 4 \left\| \tilde{u} \right\|_{L^{\infty}(p_1,p_2)} + \left\| \tilde{u}' \right\|_{L^1(p_1,p_2)} \right) m^{-1} \; .
		\end{equation}
		\item \textit{Case $\frac{x}{t} \notin f'\big(\tilde{I}\big)$.} Firstly, let us suppose $\frac{x}{t} \geqslant f'(\tilde{p}_2)$. Here there is no stationary points inside the integration interval and so it is possible to bound $\psi'$ from below by a non-zero constant,
		\begin{equation*}
			\forall \, p \in [p_1,p_2] \qquad \psi'(p) = \frac{x}{t} - f'(p) \geqslant f'(\tilde{p}_2) - f'(p_2) =: m_{\tilde{p}_2} > 0 \; .
		\end{equation*}
		Theorem \ref{3-VDC4} is then applicable and provides
		\begin{equation*}
			\forall \, t > 0 \qquad \forall \, x \geqslant f'(\tilde{p}_2) \, t \qquad \big| u_f(t,x) \big| \leqslant c_{x/t \geqslant f'(\tilde{p}_2)}(u_0,f) \, t^{-\mu} \; ,
		\end{equation*}
		with $\displaystyle c_{x/t \geqslant f'(\tilde{p}_2)}(u_0,f) := \frac{1}{2 \pi} \, \frac{1}{\mu} \, \left\| \tilde{u} \right\|_{L^{\infty}(p_1,p_2)} + \frac{1}{2 \pi} \left( 4 \left\| \tilde{u} \right\|_{L^{\infty}(p_1,p_2)} + \left\| \tilde{u}' \right\|_{L^1(p_1,p_2)} \right) m_{\tilde{p}_2}^{\, -1}$.\\ In the other case $\frac{x}{t} \leqslant f'(\tilde{p}_1)$, similar arguments furnish
		\begin{equation*}
			\forall \, t > 0 \qquad \forall \, x \leqslant f'(\tilde{p}_1) \, t \qquad \big| u_f(t,x) \big| \leqslant c_{x/t \leqslant f'(\tilde{p}_1)}(u_0,f) \, t^{-\mu} \; ,
		\end{equation*}
		with $\displaystyle c_{x/t \leqslant f'(\tilde{p}_2)}(u_0,f) := \frac{1}{2 \pi} \, \frac{1}{\mu} \, \left\| \tilde{u} \right\|_{L^{\infty}(p_1,p_2)} + \frac{1}{2 \pi} \left( 4 \left\| \tilde{u} \right\|_{L^{\infty}(p_1,p_2)} + \left\| \tilde{u}' \right\|_{L^1(p_1,p_2)} \right) m_{\tilde{p}_1}^{\; -1}$, where we set $m_{\tilde{p}_1} := f'(p_1) - f'(\tilde{p}_1) > 0$.\\ So we can finally write
		\begin{equation*}
			\forall \, (t,x) \in \mathfrak{C}_f(\tilde{p}_1,\tilde{p}_2)^c \qquad \big| u_f(t,x) \big| \leqslant c_{\tilde{I}}(u_0,f) \, t^{-\mu} \; ,
		\end{equation*}
		where
		\begin{equation} \label{c_I^c}
			c_{\tilde{I}}(u_0,f) := c_{x/t \geqslant f'(\tilde{p}_2)}(u_0,f) + c_{x/t \leqslant f'(\tilde{p}_1)}(u_0,f) \; .
		\end{equation}
	\end{enumerate}
	To prove the optimality of the above rates, we recall that the regular factor of $\tf u_0$ is supposed to be continuously differentiable on $[p_1,p_2]$. Hence Theorem \ref{3-OPTI2} is applicable and it furnishes the optimality of the rate $t^{-\mu}$ in the region $\mathfrak{C}_f(\tilde{p}_1,\tilde{p}_2)^c$. Moreover the definition of the function $\tilde{\psi}$ (see \eqref{regular-phase}) implies that this function is right continuously differentiable at $p_0$, so we can employ Theorem \ref{3-OPTI1} to prove the optimality of the rate $t^{-\frac{\mu}{2}}$ in $\mathfrak{C}_f(\tilde{p}_1,\tilde{p}_2)$. In particular, the decay rate is attained on the space-time direction defined by $\frac{x}{t} = f'(\tilde{p}_1)$.
\end{proof}

\vspace{4mm}

An $L^{\infty}$-norm estimate for the solution can be easily derived from the preceding result.

\begin{4-SCHRO2} \label{4-SCHRO2}
	Suppose that $u_0$ satisfies Condition \emph{(C1$_{[p_1,p_2],\mu}$)}. Then we have
		\begin{equation*} \label{schro21}
			\forall \, t > 0 \qquad \big\| u_f(t,.) \big\|_{L^{\infty}(\R)} \leqslant c(u_0,f) \, t^{- \frac{\mu}{2}} + c_{\tilde{I}}(u_0,f) \, t^{- \mu} \; ,
		\end{equation*}
		where the constants $c(u_0,f) > 0$ and $c_{\tilde{I}}(u_0,f) > 0$ are given by \eqref{c_I} and \eqref{c_I^c} respectively. In particular, we have
		\begin{equation*} \label{schro22}
			\forall \, t \geqslant 1 \qquad \big\| u_f(t,.) \big\|_{L^{\infty}(\R)} \leqslant \big( c(u_0,f) + c_{\tilde{I}}(u_0,f) \big) \, t^{- \frac{\mu}{2}} \; .
		\end{equation*}
\end{4-SCHRO2}

\begin{proof}
	Simple consequence of Theorem \ref{4-SCHRO1}.
\end{proof}

\vspace{4mm}

As an application of the above theorem, we provide an $L^{\infty}$-norm estimate of the solution of the free Schrödinger equation on the line for initial data satisfying Condition (C1$_{[p_1,p_2],\mu}$). The resulting decay rate is given by $t^{-\frac{\mu}{2}}$. Let us remark that this decay rate has been obtained in our result \cite[Theorem 5.6]{fusion} by expanding the solution to one term on the space-time direction given by $\frac{x}{t} = 2 \, p_1$.

\begin{4-OPTI} \label{4-OPTI}
	Let $u_S : \R_+ \times \R \longrightarrow \C$ be the solution of the free Schrödinger equation on $\R$,
	\begin{equation*}
	\left\{ \begin{array}{rl}
			& \hspace{-2mm} \big[ i \, \partial_t + \partial_{xx} \big] u_S(t) = 0 \\ [2mm]
			& \hspace{-2mm} u_S(0) = u_0
	\end{array} \right. \; ,
\end{equation*}
	for $t \geqslant 0$, where $u_0$ satisfies Condition \emph{(C1$_{[p_1,p_2],\mu}$)}. Then we have
	\begin{equation*} \label{optischro}
		\forall \, t \geqslant 1 \qquad \big\| u_S(t,.) \big\|_{L^{\infty}(\R)} \leqslant c(u_0,f_S) \, t^{- \frac{\mu}{2}} \; ,
	\end{equation*}
	where the constant $c(u_0,f_S) > 0$ can be computed from Theorem \ref{4-SCHRO2}, and the decay rate is optimal.
\end{4-OPTI}

\begin{proof}
	Application of Theorem \ref{4-SCHRO2} for the symbol $f_S(p) = p^2$, which gives the differential operator $- \partial_{xx}$.
\end{proof}

\vspace{4mm}

The aim of the two following theorems is to show that the influence of the singular frequency $p_1$ on the decay rate is stronger in space-time regions containing the space-time direction $\frac{x}{t} = f'(p_1)$.\\
In the following result, we establish estimates of the solution in arbitrarily narrow cones containing the above space-time direction. In such regions, the phase function, coming from the rewriting of the solution as an oscillatory integral, has a stationary point which is in a neighbourhood of the singular frequency $p_1$. In this context, these two particular points are expected to interact with each other, producing the slow decay $t^{-\frac{\mu}{2}}$.\\
Here we do not require the initial data to be in a frequency band anymore. This permits to remove the localization phenomenon produced by the frequency band, which has been exhibited in Theorem \ref{4-SCHRO1}, and to focus only on the influence of the singular frequency $p_1$ on the decay rate in the above mentioned cones.\\

\noindent \textbf{Condition (C2$_{p_1, \mu}$).} Fix $\mu \in (0,1)$ and choose a finite real number $p_1$.\\ A tempered distribution $u_0$ on $\R$ satisfies Condition (C2$_{p_1, \mu}$) if and only if $\tf u_0 \in L^1(\R)$ and there exists a bounded differentiable function $\tilde{u} : \R \longrightarrow \C$ such that $\tilde{u}(p_1) \neq 0$, $\tilde{u}' \in L^1(\R)$ and
\begin{equation*}
	\forall \, p \in \R \backslash \{ p_1 \} \qquad \tf u_0(p) = |p-p_1|^{\mu-1} \, \tilde{u}(p) \; .
\end{equation*}

\vspace{2mm}

\begin{4-REM9} \em
	\begin{enumerate}
		\item One can follow the lines of the point i) of Remark \ref{4-REM8} to ensure that the subset of tempered distributions satisfying Condition (C2$_{p_1,\mu}$) is non-empty.
		\item Here $u_0$ is at least a continuous function on $\R$ but it is not necessarily analytic. Furthermore the solution formula \eqref{formula2} is still well-defined for all $t \geqslant 0$ and $x \in \R$.
	\end{enumerate}
\end{4-REM9}

\vspace{1mm}

\begin{4-PREGEO1} \label{4-PREGEO1}
	Suppose that $u_0$ satisfies Condition \emph{(C2$_{p_1, \mu}$)} and choose a finite real number $\varepsilon > 0$. Then for all $\displaystyle (t,x) \in \mathfrak{C}_f(p_1 - \varepsilon, p_1 + \varepsilon)$, we have
		\begin{equation*}
			\big| u_f(t,x) \big| \leqslant c^{(1)}(u_0,f) \, t^{-\frac{\mu}{2}} + c_{ \varepsilon}^{(2)}(u_0,f) \, t^{-1} \; .
		\end{equation*}
		The constants $c^{(1)}(u_0,f)$ and $c_{\varepsilon}^{(2)}(u_0,f)$ are given by \eqref{c_1} and \eqref{c_2} respectively.
\end{4-PREGEO1}

\begin{proof}
	We shall employ the rewritting of the solution given in the proof of Theorem \ref{4-SCHRO1}, i.e.
	\begin{equation*}
		\forall \, (t,x) \in (0,+\infty) \times \R \qquad u_f(t,x) = \int_{\R} U(p) \, e^{i t \psi(p)} \, dp \; ,
	\end{equation*}
	where
	\begin{equation*}
		\left\{ \begin{array}{rl}
				& \displaystyle \forall \, p \in \R \backslash \{ p_1 \} \qquad U(p) := \frac{1}{2\pi}  \, \tf u_0(p) = \frac{1}{2\pi}  \, |p-p_1|^{\mu-1} \, \tilde{u}(p) \; , \\
				& \vspace{-0.3cm} \\
				& \displaystyle \forall \, p \in \R \qquad \psi(p) := \frac{x}{t} \, p - f(p) \; .
		\end{array} \right.
	\end{equation*}
	Let $\varepsilon > 0$, choose a finite real number $\eta > 0$ such that $\eta > \varepsilon$ (for example $\eta = \varepsilon + 1$) and split the integral as follows,
	\begin{align*}
		\int_{\R} U(p) \, e^{i t \psi(p)} \, dp	& = \int_{p_1- \eta}^{p_1 + \eta} \dots \; + \int_{\R \backslash [p_1 - \eta, p_1 + \eta]} \dots =: I^{(1)}(t,x,\eta) + I^{(2)}(t,x,\eta) \; .
	\end{align*}
	Firstly we study $I^{(1)}(t,x,\eta)$. We recall that
	\begin{equation*}
		\psi'(p) = \frac{x}{t} - f'(p) \; ;
	\end{equation*}
	since $\displaystyle \frac{x}{t}$ is supposed to belong to $\displaystyle \big( f'(p_1 - \varepsilon), f'(p_1 + \varepsilon) \big)$, then $\psi$ has a stationary point which belongs to $(p_1-\varepsilon, p_1 + \varepsilon) \subset [p_1-\eta, p_1 + \eta]$. Following the arguments of the point i) of the proof of Theorem \ref{4-SCHRO1}, we apply Theorem \ref{3-VDC} on $[p_1 - \eta, p_1]$ and on $[p_1, p_1 + \eta]$ with $\rho = 2$, leading to
	\begin{equation*}
		\Big| I^{(1)}(t,x,\eta) \Big| \leqslant \left| \int_{p_1-\eta}^{p_1} \dots \right| + \left| \int_{p_1}^{p_1 + \eta} \dots \right| \leqslant \left( c_{1}^{(1)}(u_0,f) + c_{2}^{(1)}(u_0,f) \right) t^{-\frac{\mu}{2}} \; ,
	\end{equation*}
	where
	\begin{align*}
		& \bullet \quad c_{1}^{(1)}(u_0,f) := \frac{1}{2 \pi} \, \frac{3}{\mu} \, \left\| \tilde{u} \right\|_{L^{\infty}(p_1-\eta, p_1)} + \frac{1}{\pi} \left( 4 \left\| \tilde{u} \right\|_{L^{\infty}(p_1-\eta, p_1)} + \left\| \tilde{u}' \right\|_{L^1(p_1-\eta, p_1)} \right) m_{1,\eta}^{\; -1} \; , \\
		& \bullet \quad c_{2}^{(1)}(u_0,f) := \frac{1}{2 \pi} \, \frac{3}{\mu} \, \left\| \tilde{u} \right\|_{L^{\infty}(p_1, p_1 + \eta)} + \frac{1}{\pi} \left( 4 \left\| \tilde{u} \right\|_{L^{\infty}(p_1, p_1 + \eta)} + \left\| \tilde{u}' \right\|_{L^1(p_1, p_1 + \eta)} \right) m_{2,\eta}^{\; -1} \; ,
	\end{align*}
	with $\displaystyle m_{1, \eta} := \min_{p \in [p_1-\eta, p_1]} f''(p) > 0$ and $\displaystyle m_{2, \eta} := \min_{p \in [p_1, p_1 + \eta]} f''(p) > 0$. The constant $c^{(1)}(u_0, f)$ is then defined by
	\begin{equation} \label{c_1}
		c^{(1)}(u_0, f) := c_{1}^{(1)}(u_0,f) + c_{2}^{(1)}(u_0,f) \; .
	\end{equation}
	Let us study $I^{(2)}(t,x,\eta)$. Let $k \in \N$ and consider the following sequence,
	\begin{equation*}
		\tilde{I}_k^{(2)}(t,x, \eta) := \int_{p_1+\eta}^{p_1+ \eta + k} U(p) \, e^{it \psi(p)} \, dp \; .
	\end{equation*}
	Since $\displaystyle \frac{x}{t} \in \big( f'(p_1 - \varepsilon) , f'(p_1 + \varepsilon)  \big)$, we note that the first derivative of the phase function does not vanish on $[p_1 + \eta, p_1 + \eta + k]$ and more precisely, we have for any $k \in \N$,
	\begin{equation*}
		\forall \, p \in [p_1 + \eta, p_1 + \eta + k] \quad \big| \psi'(p) \big| = f'(p) - \frac{x}{t} \geqslant f'(p_1 + \eta) - f'(p_1 + \varepsilon) =: \tilde{m}_{1, \eta, \varepsilon} > 0  \; .
	\end{equation*}
	Theorem \ref{3-VDC4} in the case $\mu = 1$ furnishes for all $(t,x) \in \mathfrak{C}_f(p_1 - \varepsilon, p_1 + \varepsilon)$,
	\begin{equation*}
		\left| \tilde{I}_k^{(2)}(t,x, \eta) \right| \leqslant \frac{1}{2\pi} \left( 3 \, \left\| U \right\|_{L^{\infty}(p_1+\eta, p_1 + \eta + k)} +  \left\| U' \right\|_{L^1(p_1+\eta, p_1 + \eta + k)} \right) \tilde{m}_{1, \eta, \varepsilon}^{\; -1} \, t^{-1} \; .
	\end{equation*}
	Now by using the hypotheses on the initial data, we give estimates for $\left\| U \right\|_{L^{\infty}(p_1+\eta, p_1 + \eta + k)}$ and $\left\| U' \right\|_{L^1(p_1+\eta, p_1 + \eta + k)}$, namely,
	\begin{align*}
		& \bullet \quad \forall \, p \in [p_1 + \eta, p_1 + \eta + k] \qquad \left| U(p) \right| \leqslant \eta^{\mu-1} \left\| \tilde{u} \right\|_{L^{\infty}(\R)} \; , \\[2mm]
		& \bullet \quad \int_{p_1 + \eta}^{p_1 + \eta + k} \big| U'(p) \big| dp \leqslant \eta^{\mu-1} \left( \| \tilde{u} \|_{L^{\infty}(\R)} + \| \tilde{u}' \|_{L^1(\R)} \right) \; .
	\end{align*}
	Consequently, $\tilde{I}_2^{(k)}(t,x,\eta)$ can be estimated as follows,
	\begin{equation} \label{tildeI2}
		\left| \tilde{I}_k^{(2)}(t,x, \eta) \right| \leqslant \frac{1}{2\pi} \, \eta^{\mu-1} \left( 4 \, \left\| \tilde{u} \right\|_{L^{\infty}(\R)} +  \left\| \tilde{u}' \right\|_{L^1(\R)} \right) \tilde{m}_{1, \eta, \varepsilon}^{\; -1} \, t^{-1} \; .			
	\end{equation}
	Using the dominated convergence Theorem which claims that
	\begin{equation*}
		\lim_{k \rightarrow + \infty} \tilde{I}_k^{(2)}(t,x,\eta) = \int_{p_1 + \eta}^{+\infty} U(p) \, e^{i t \psi(p)} dp \; ,
	\end{equation*}
	we can take the limit in \eqref{tildeI2} providing
	\begin{equation*}
		\left| \int_{p_1 + \eta}^{+ \infty} U(p) \, e^{i t \psi(p)} dp \right| \leqslant c_{1,\varepsilon}^{(2)}(u_0, f) \, t^{-1} \; ,
	\end{equation*}
	with
	\begin{equation*}
		c_{1,\varepsilon}^{(2)}(u_0, f) := \frac{1}{2\pi} \, \eta^{\mu-1} \left( 4 \, \left\| \tilde{u} \right\|_{L^{\infty}(\R)} +  \left\| \tilde{u}' \right\|_{L^1(\R)} \right) \tilde{m}_{1, \eta, \varepsilon}^{\; -1} \; .
	\end{equation*}
	Similar arguments furnish the following estimate,
	\begin{equation*}
		\forall \, (t,x) \in \mathfrak{C}_f (p_1-\varepsilon, p_1 + \varepsilon) \qquad \left| \int_{- \infty}^{p_1 - \eta} U(p) \, e^{i t \psi(p)} dp \right| \leqslant c_{2,\varepsilon}^{(2)}(u_0, f) \, t^{-1} \; ,
	\end{equation*}
	with
	\begin{equation*}
		c_{2,\varepsilon}^{(2)}(u_0, f) := \frac{1}{2\pi} \, \eta^{\mu-1} \left( 4 \, \left\| \tilde{u} \right\|_{L^{\infty}(\R)} +  \left\| \tilde{u}' \right\|_{L^1(\R)} \right) \tilde{m}_{2, \eta, \varepsilon}^{\; -1} \; ,
	\end{equation*}
	where $\displaystyle \tilde{m}_{2, \eta, \varepsilon} := f'(p_1 - \varepsilon) - f'(p_1 - \eta) > 0$.\\
	Finally, by setting
	\begin{equation} \label{c_2}
		c_{\varepsilon}^{(2)}(u_0, f) := c_{1,\varepsilon}^{(2)}(u_0, f) + c_{2,\varepsilon}^{(2)}(u_0, f) \; ,
	\end{equation}
	we obtain for all $(t,x) \in \mathfrak{C}_f(p_1-\varepsilon, p_1 + \varepsilon)$,
	\begin{equation*}
		\Big| I^{(2)}(t,x,\eta) \Big| \leqslant \left| \int_{-\infty}^{p_1-\eta} \dots \right| + \left| \int_{p_1 + \eta}^{+\infty} \dots \right| \leqslant c_{\varepsilon}^{(2)}(u_0, f) \, t^{-1} \; ,
	\end{equation*}
	which ends the proof.
\end{proof}

\vspace{4mm}

Now we provide estimates of the solution in space-time cones which do not contain the critical direction given by the singular frequency. In this case, the distance between the stationary point and the singular frequency is bounded from below by a positive constant, which removes the superposition of the effects of these particular points. Hence these two points provide two distinct decay rates: $t^{-\frac{1}{2}}$ coming from the stationary point and $t^{-\mu}$ coming from the singular frequency. We note that these two rates are better than $t^{-\frac{\mu}{2}}$.\\
Theorems \ref{4-PREGEO2} and \ref{4-PREGEO1} can be compared with Theorem 5.7 in \cite{fusion}. In the latter, we have furnished estimates of the solution of the free Schrödinger equation in space-time regions along the direction $\frac{x}{t} = 2 \, p_1$, and this direction is outside the regions. The estimates show that the decay rate diminishes when the boundary of the region approaches the direction given by $p_1$.

\begin{4-PREGEO2} \label{4-PREGEO2}
	Suppose that $u_0$ satisfies Condition \emph{(C2$_{p_1, \mu}$)} and choose two finite real numbers $\tilde{p}_1 < \tilde{p}_2$ such that $p_1 \notin [\tilde{p}_1, \tilde{p}_2]$. Then for all $\displaystyle (t,x) \in \mathfrak{C}_f(\tilde{p}_1, \tilde{p}_2)$, we have
		\begin{equation*}
			\big| u_f(t,x) \big| \leqslant c_{\tilde{p}_1, \tilde{p}_2}^{(1)}(u_0,f) \, t^{-\frac{1}{2}} + c_{\tilde{p}_1, \tilde{p}_2}^{(2)}(u_0,f) \, t^{-\mu} + c_{\tilde{p}_1, \tilde{p}_2}^{(3)}(u_0,f) \, t^{-1} \; .
		\end{equation*}
		The constants $c_{\tilde{p}_1, \tilde{p}_2}^{(1)}(u_0,f)$, $c_{\tilde{p}_1, \tilde{p}_2}^{(2)}(u_0,f)$ and $c_{\tilde{p}_1, \tilde{p}_2}^{(3)}(u_0,f)$ are given by \eqref{c_21}, \eqref{c_22} and \eqref{c_23} respectively.
\end{4-PREGEO2}

\begin{proof}
	The calculations in the present proof are similar to those of the proof of Theorem \ref{4-PREGEO1}. Thus we give only the main steps of the proof of Theorem \ref{4-PREGEO2}. \\
	Let $\eta \in \big(0, \min \{ |\tilde{p}_1 - p_1|, |\tilde{p}_2 - p_1| \} \big)$ and split the integral as follows,
	\begin{align*}
		u_f(t,x) = \int_{\R} U(p) \, e^{it \psi(p)} \, dp	& = \int_{\tilde{p}_1 - \eta}^{\tilde{p}_2 + \eta} \dots \; + \int_{\R \backslash [\tilde{p}_1 - \eta, \tilde{p}_2 + \eta]} \dots \\[3mm]
												& =: I^{(1)}(t,x,\eta) + I^{(2)}(t,x, \eta) \; ,
	\end{align*}
	On the interval $[\tilde{p}_1- \eta, \tilde{p}_2 + \eta]$, the phase has a unique stationary point and the amplitude has no singular points. Theorem \ref{3-VDC} is applicable with $\rho = 2$ and $\mu = 1$, and it leads to
	\begin{equation*}
		\forall \, (t,x) \in \mathfrak{C}_f(\tilde{p}_1, \tilde{p}_2) \qquad \Big| I^{(1)}(t,x,\eta) \Big| \leqslant c_{\tilde{p}_1, \tilde{p}_2}^{(1)}(u_0,f) \, t^{-\frac{1}{2}} \; ,
	\end{equation*}
	where
	\begin{equation} \label{c_21}
		c_{\tilde{p}_1, \tilde{p}_2}^{(1)}(u_0,f) := \left\{ \begin{array}{rl}
				& \displaystyle \frac{(\tilde{p}_1 - \eta - p_1)^{\mu-1}}{\pi} \, \bigg( \left\| \tilde{u} \right\|_{L^{\infty}(\R)} \\
				& \displaystyle \qquad + \; \Big( 4  \left\| \tilde{u} \right\|_{L^{\infty}(\R)} + \left\| \tilde{u}' \right\|_{L^1(\R)} \Big) m_{1, \tilde{p}_1, \tilde{p}_2}^{\; -1} \bigg) \; , \quad \text{if} \; p_1 < \tilde{p}_1 \; , \\
				& \vspace{-0.3cm} \\
				& \displaystyle \frac{(p_1 - \tilde{p}_2 - \eta)^{\mu-1}}{\pi} \, \bigg( \left\| \tilde{u} \right\|_{L^{\infty}(\R)} \\
				& \displaystyle \qquad + \; \Big( 4  \left\| \tilde{u} \right\|_{L^{\infty}(\R)} + \left\| \tilde{u}' \right\|_{L^1(\R)} \Big) m_{1, \tilde{p}_1, \tilde{p}_2}^{\; -1} \bigg) \; , \quad \text{if} \; p_1 > \tilde{p}_2 \; ,
		\end{array} \right.
	\end{equation}
	with $\displaystyle m_{1, \tilde{p}_1, \tilde{p}_2} := \min_{p \in [\tilde{p}_1 - \eta, \tilde{p}_2 + \eta]} f''(p) > 0$.\\
	Now let us study $I^{(2)}(t,x,\eta)$. First of all, we remark that we integrate over two infinite intervals such that one of them contains the singular frequency $p_1$. Consequently we shall suppose that $p_1 < \tilde{p}_1$ without loss of generality; the other case $p_1 > \tilde{p}_2$ can be treated in a similar way.\\
	We consider the following sequence,
	\begin{equation*}
		\forall \, k \in \N^* \qquad \tilde{I}_k^{(2)}(t,x,\eta) := \int_{p_1-k}^{\tilde{p}_1-\eta} U(p) \, e^{it \psi(p)} \, dp \; .
	\end{equation*}
	We note that $[p_1 - k, \tilde{p}_1 - \eta]$ contains the singular frequency $p_1$ and $\psi'$ does not vanish on this interval, and thus Theorem \ref{3-VDC4} is applicable on $[p_1 - k, p_1]$ and on $[p_1, \tilde{p}_1 - \eta]$. Then we take the limit when $k$ tends to infinity by using the dominated convergence Theorem and we obtain
	\begin{equation} \label{I_k-2}
		\forall \, (t,x) \in \mathfrak{C}_f(\tilde{p}_1, \tilde{p}_2) \qquad \left| \int_{-\infty}^{\tilde{p}_1-\eta} U(p) \, e^{it \psi(p)} \, dp \right| \leqslant c_{\tilde{p}_1, \tilde{p}_2}^{(2)}(u_0,f) \, t^{-\mu} \; ,
	\end{equation}
	where
	\begin{equation*}
		c_{\tilde{p}_1, \tilde{p}_2}^{(2)}(u_0,f) := \frac{1}{\pi} \left( \frac{1}{\mu} \left\| \tilde{u} \right\|_{L^{\infty}(\R)} + \Big( 4 \left\| \tilde{u} \right\|_{L^{\infty}(\R)} + \left\| \tilde{u}' \right\|_{L^1(\R)} \Big) \tilde{m}_{2, \tilde{p}_1}^{\; -1} \right) \; ,
	\end{equation*}	
	with $\displaystyle \tilde{m}_{2, \tilde{p}_1} := f'(\tilde{p}_1) - f'(\tilde{p}_1 - \eta) > 0$. To study the integral on the other infinite interval, we define $\tilde{I}_k^{(3)}(t,x,\eta)$ as follows,
	\begin{equation*}
		\forall \, k \in \N^* \qquad \tilde{I}_k^{(3)}(t,x,\eta) := \int_{\tilde{p}_2+\eta}^{\tilde{p}_2+\eta+k} U(p) \, e^{it \psi(p)} \, dp \; .
	\end{equation*}
	Here there is no singular frequency or stationary point, therefore Theorem \ref{3-VDC4} in the case $\mu = 1$ is applicable and it furnishes
	\begin{equation} \label{I_k-3}
		\forall \, (t,x) \in \mathfrak{C}_f(\tilde{p}_1), \tilde{p}_2) \qquad \left| \int_{\tilde{p}_2+\eta}^{+\infty} U(p) \, e^{it \psi(p)} \, dp \right| \leqslant c_{\tilde{p}_1, \tilde{p}_2}^{(3)}(u_0,f) \, t^{-1} \; ,
	\end{equation}
	where
	\begin{equation*}
		c_{\tilde{p}_1, \tilde{p}_2}^{(3)}(u_0,f) := \frac{(\tilde{p}_2 + \eta - p_1)^{\mu-1}}{2 \pi} \left( 4 \left\| \tilde{u} \right\|_{L^{\infty}(\R)} + \left\| \tilde{u}' \right\|_{L^1(\R)} \right) \tilde{m}_{3, \tilde{p}_2}^{\; -1} \; ,
	\end{equation*}
	with $\displaystyle \tilde{m}_{3, \tilde{p}_2} := f'(\tilde{p}_2 + \eta) - f'(\tilde{p}_2) > 0$.\\
	Consequently we derive the following estimate for $I^{(2)}(t,x,\eta)$ in the case $p_1 < \tilde{p}_1$ from \eqref{I_k-2} and \eqref{I_k-3},
	\begin{equation*}
		\forall \, (t,x) \in \mathfrak{C}_f(\tilde{p}_1, \tilde{p}_2) \qquad \Big| I^{(2)}(t,x,\eta) \Big| \leqslant c_{\tilde{p}_1, \tilde{p}_2}^{(2)}(u_0,f) \, t^{-\mu} + c_{\tilde{p}_1, \tilde{p}_2}^{(3)}(u_0,f) \, t^{-1} \; .
	\end{equation*}
	To conclude, we provide the values of the constants $c_{\tilde{p}_1, \tilde{p}_2}^{(2)}(u_0,f)$ and $c_{\tilde{p}_1, \tilde{p}_2}^{(3)}(u_0,f)$ depending on the position of $p_1$ with respect to the interval $[\tilde{p}_1, \tilde{p}_2]$:
	\begin{equation} \label{c_22}
		\bullet \quad c_{\tilde{p}_1, \tilde{p}_2}^{(2)}(u_0,f) := \left\{ \begin{array}{rl}
				& \displaystyle \frac{1}{\pi} \left( \frac{1}{\mu} \left\| \tilde{u} \right\|_{L^{\infty}(\R)} + \Big( 4 \left\| \tilde{u} \right\|_{L^{\infty}(\R)} + \left\| \tilde{u}' \right\|_{L^1(\R)} \Big) \tilde{m}_{2, \tilde{p}_1}^{\; -1} \right) \; , \quad \text{if} \; p_1 < \tilde{p}_1 \; , \\
				& \vspace{-0.3cm} \\
				& \displaystyle \frac{1}{\pi} \left( \frac{1}{\mu} \left\| \tilde{u} \right\|_{L^{\infty}(\R)} + \Big( 4 \left\| \tilde{u} \right\|_{L^{\infty}(\R)} + \left\| \tilde{u}' \right\|_{L^1(\R)} \Big) \tilde{m}_{3, \tilde{p}_2}^{\; -1} \right) \; , \quad \text{if} \; p_1 > \tilde{p}_2 \; ,
		\end{array} \right.
	\end{equation}
	\begin{equation} \label{c_23}
		\bullet \qquad c_{\tilde{p}_1, \tilde{p}_2}^{(3)}(u_0,f) := \left\{ \begin{array}{rl}
				& \displaystyle \frac{(\tilde{p}_2 + \eta - p_1)^{\mu-1}}{2 \pi} \left( 4 \left\| \tilde{u} \right\|_{L^{\infty}(\R)} + \left\| \tilde{u}' \right\|_{L^1(\R)} \right) \tilde{m}_{3, \tilde{p}_2}^{\; -1} \; , \quad \text{if} \; p_1 < \tilde{p}_1 \; , \\
				& \vspace{-0.3cm} \\
				& \displaystyle \frac{(p_1 - \tilde{p}_1 + \eta)^{\mu-1}}{2 \pi} \left( 4 \left\| \tilde{u} \right\|_{L^{\infty}(\R)} + \left\| \tilde{u}' \right\|_{L^1(\R)} \right) \tilde{m}_{2, \tilde{p}_1}^{\; -1} \; , \quad \text{if} \; p_1 > \tilde{p}_2 \; .
		\end{array} \right.
	\end{equation}
\end{proof}

\section{An intrinsic localization phenomenon caused by a limited growth of the symbol} \label{Sec-App2}

\hspace{0.5cm} As explained in Section \ref{Sec-results}, a symbol having a first derivative which is bounded may influence the dispersion of the solution of equation \eqref{evoleq}: if $f'(\R) = (a,b)$, where $a < b$ are two finite real numbers given by the limits of $f'$ at $-\infty$ and $+\infty$, and if the initial datum is in a bounded frequency band $[p_1,p_2]$, then the associated solution will have a limited speed when the time tends to infinity. This is a consequence of the time-asymptotic localization of the solution in $\mathfrak{C}_f(p_1,p_2)$ and of the following inclusion coming from the boundedness of $f'$:
\begin{equation*}
	\mathfrak{C}_f(p_1,p_2) \subset \left\{ (t,x) \in (0,+\infty) \times \R \, \Big| \, a < \frac{x}{t} < b \right\} =: \mathfrak{C}_f(-\infty,+\infty) \; ,
\end{equation*}
This can be viewed as an asymptotic version of the notion of causality for initial data without compact support. Our aim in this section consists in extending this phenomenon to the case of initial data which are not in bounded frequency bands; this is stated in Theorem \ref{5-SCHRO3}.  An illustration will be given in the setting of the Klein-Gordon equation on the line in Corollary \ref{5-COR3}.\\

We recall that we study the class of evolution equations on the line given by
\begin{equation*}
	\left\{ \begin{array}{rl}
			& \hspace{-2mm} \left[ i \, \partial_t - f \big(D\big) \right] u_f(t) = 0 \\ [2mm]
			& \hspace{-2mm} u_f(0) = u_0
	\end{array} \right. \; ,
\end{equation*}
for $t \geqslant 0$. In the present section, we make technical hypotheses on the behaviour of $f''$ at infinity instead of assuming only the boundedness of $f$': this is required by the method we use to extend the above intrinsic localization. However these hypotheses imply in particular that $f'(\R) = (a,b)$, where $a$ and $b$ are the limits of $f'$ at infinity, and that the distance between $f'$ and its limits at infinity can be estimated from below; see Lemma \ref{5-LEM3}:\\

\noindent \textbf{Condition (S$_{\beta_+, \beta_-,R}$).} Fix $\beta_- \geqslant \beta_+ > 1$ and $R \geqslant 1$.\\
A $\mathcal{C}^{\infty}$-function $f : \R \longrightarrow \R$ satisfies Condition (S$_{\beta_+, \beta_-,R}$) if and only if the second derivative of $f$ verifies $f'' > 0$ and
\begin{equation} \label{condS}
	\exists \, c_+ \geqslant c_- > 0 \qquad \forall  \, |p| \geqslant R \qquad c_- \, |p|^{-\beta_-} \leqslant f''(p) \leqslant c_+ \, |p|^{-\beta_+} \; .
\end{equation}

\vspace{2mm}

\begin{5-LEM3} \label{5-LEM3}
	Let $f : \R \longrightarrow \R$ be a function satisfying Condition \emph{(S$_{\beta_+, \beta_-,R}$)}. Then
	\begin{enumerate}
		\item we have $f'(\R) = (a,b)$ where
		\begin{equation*}
			a := \lim_{p \rightarrow - \infty} f'(p) \qquad , \qquad b := \lim_{p \rightarrow + \infty} f'(p) \; ;
		\end{equation*}
		\item we have
		\begin{align*}
			& \bullet \quad \forall \, p \geqslant R \qquad b - f'(p) \geqslant \frac{c_-}{\beta_- - 1} \, p^{1-\beta_-} \; , \\
			& \bullet \quad \forall \, p \leqslant -R \qquad f'(p) - a \geqslant \frac{c_-}{\beta_- - 1} \, (-p)^{1-\beta_-} \; .
		\end{align*}
	\end{enumerate}
\end{5-LEM3}

\begin{proof}
	\begin{enumerate}
		\item On the compact interval $[-R,R]$, the function $f'$ is bounded since it is continuous. Now, using the right inequality in \eqref{condS}, we have for $p \geqslant R$,
		\begin{equation*}
			f'(p) - f'(R) = \int_R^p f''(x) \, dx \leqslant c_+ \int_R^p x^{-\beta_+} dx \; ,
		\end{equation*}
		which provides
		\begin{equation*}
		f'(p) \leqslant \frac{c_+}{1-\beta_+} \, p^{1-\beta_+} + f'(R) - \frac{c_+}{1 - \beta_+} \, R^{1-\beta_+} \leqslant f'(R) - \frac{c_+}{1 - \beta_+} \, R^{1-\beta_+} < \infty \; .
		\end{equation*}
		Consequently $f'$ is bounded from above on $[R, +\infty)$ and similar arguments show that $f'$ is bounded from below on $(-\infty, R]$. Since the function $f'$ is strictly increasing on $\R$, we deduce that $f'$ is bounded on $\R$ and its bounds are given by its limits at $-\infty$ and $+\infty$.		
		\item For $p \geqslant R$, we have
		\begin{equation*}
			b - f'(p) = \int_p^{+\infty} f''(x) \, dx \geqslant c_- \int_p^{+\infty} x^{-\beta_-} dx = - \frac{c_-}{1 - \beta_-} \, p^{1-\beta_-} \; ,
		\end{equation*}
		where we used the left inequality of \eqref{condS}. In the same way, we have for all $p \leqslant -R$,
		\begin{equation*}
			f'(p) - a = \int_{-\infty}^p f''(x) \, dx \geqslant c_- \int_{-\infty}^p (-x)^{-\beta_-} dx = - \frac{c_-}{1 - \beta_-} \, (-p)^{1-\beta_-} \; .
		\end{equation*}
	\end{enumerate}
\end{proof}

\vspace{4mm}

In Theorem \ref{5-SCHRO3}, we do not assume that the initial datum is in a bounded frequency band but we allow it to have a singular frequency, as in Theorem \ref{4-PREGEO1} and Theorem \ref{4-PREGEO2}. In this case, we put the singular frequency at $0$ in order to avoid a proof with too many technical calculations. Without this assumption on the position of the singular frequency, the result of Theorem \ref{5-SCHRO3} remains unchanged and its proof follows the steps of the proof in the case of the singular frequency at $0$. In addition to this, we make precise hypotheses on the decay of $\tf u_0$ at infinity to carry out the proof.\\

\noindent \textbf{Condition (C3$_{\mu, \alpha,r}$).} Fix $\mu \in (0,1]$, $\alpha > \mu$ and $r \geqslant 0$.\\
A tempered distribution $u_0$ on $\R$ satisfies Condition (C3$_{\mu, \alpha,r}$) if and only if there exists a bounded differentiable function $\tilde{u} : \R \longrightarrow \C$ such that $\tilde{u}(0) \neq 0$ if $\mu \neq 1$, with
\begin{equation*}
	\forall \, p \in \R \backslash \{ 0 \} \qquad \tf u_0(p) = |p|^{\mu-1} \, \tilde{u}(p) \; .
\end{equation*}
Moreover we suppose that
\begin{equation*}
	\exists \,  M \geqslant 0 \quad \forall \, p \in \R \qquad \big| \tilde{u}(p) \big| \leqslant M \left( 1 + p^2 \right)^{-\frac{\alpha}{2}} \; ,
\end{equation*}
and that $\displaystyle \tilde{u}' \in L_{loc}^1(\R)$ with
\begin{equation*}
	\exists \, M' \geqslant 0 \quad \forall \, n \in \left\{ n \in \mathbb{Z} \, \big| \, |n| \geqslant r \right\} \qquad \left\| \tilde{u}' \right\|_{L^1(n,n+1)} \leqslant M' \, |n|^{-\alpha} \; .
\end{equation*}

\vspace{1mm}

\begin{5-REM4}
	\em
	\begin{enumerate}
		\item The above condition implies in particular that $\tf u_0$ belongs to $L^1(\R)$. Indeed $\displaystyle \tf u_0 \in L_{loc}^1(\R)$ since the function $p \longmapsto |p|^{\mu-1}$ with $\mu \in (0,1]$ belongs to $L_{loc}^1(\R)$, and $\displaystyle \tilde{u} \in L^{\infty}(\R)$. Furthermore we have
	\begin{equation*}
		\forall \, p \in \R \backslash \{0 \} \qquad \big| \tf u_0(p) \big| \leqslant M \left( 1 + p^2 \right)^{-\frac{\alpha}{2}} |p|^{\mu-1} \leqslant M \, |p|^{\mu-1-\alpha} \; ,
	\end{equation*}
	since $\big( 1 + p^2 \big)^{\frac{1}{2}} \geqslant |p|$. Hence the hypothesis $\alpha > \mu$ leads to the integrability of $\tf u_0$ on $\R$.\\
	Thanks to that, one can show that the subset of tempered distributions verifying Condition (C3$_{\mu, \alpha,r}$) is non-empty by following the lines of Remark \ref{4-REM8} i), and the solution formula \eqref{formula2} is well-defined for all $t \geqslant 0$ and $x \in \R$. 
	\item Let us give an example for the above condition. Choose $u_0 \in \mathcal{S}'(\R)$ such that its Fourier transform has the following form:
	\begin{equation*}
		\forall \, p \in \R \backslash \{ 0 \} \qquad \tf u_0 (p) = |p|^{\mu-1} \, (1+p^2)^{-\frac{\alpha}{2}} \; ,
	\end{equation*}
	with $\mu \in (0,1]$ and $\alpha > \mu$. Here $\tilde{u} : \R \longrightarrow \R$ is given by $\displaystyle \tilde{u}(p) = (1+p^2)^{-\frac{\alpha}{2}}$ for all $p \in \R$. \\
	In this case, we only have to control $\left\| \tilde{u}' \right\|_{L^1(n,n+1)}$ since the other hypotheses are clearly satisfied. One can quickly show that
	\begin{equation*}
		\left\| \tilde{u}' \right\|_{L^1(n,n+1)} = \left\{ \begin{array}{rl}
				& \displaystyle \tilde{u}(n) - \tilde{u}(n+1) \leqslant \tilde{u}(n) \; , \quad \text{if} \; n \geqslant 0 \; , \\
				& \vspace{-0.3cm} \\
				& \displaystyle \tilde{u}(n+1) - \tilde{u}(n) \leqslant \tilde{u}(n+1) \; , \quad \text{if} \; n \leqslant -1 \; .
		\end{array} \right.
	\end{equation*}
	Using the fact that $|n+1|^{-\alpha} \leqslant 2^{\alpha} |n|^{-\alpha}$, if $n \leqslant -2$ according to Lemma \ref{6-LEM2} (see below), we obtain
	\begin{equation*}
		\left\| \tilde{u}' \right\|_{L^1(n,n+1)} \leqslant \left\{ \begin{array}{rl}
				& \displaystyle \left( 1 + n^2 \right)^{-\frac{\alpha}{2}} \leqslant n^{-\alpha} \; , \quad \text{if} \; n \geqslant 0 \; , \\
				& \vspace{-0.3cm} \\
				& \displaystyle \left( 1 + (n+1)^2 \right)^{-\frac{\alpha}{2}} \leqslant |n+1|^{-\alpha} \leqslant 2^{\alpha} |n|^{-\alpha} \; , \quad \text{if} \; n \leqslant -2 \; .
		\end{array} \right.
	\end{equation*}
	Hence for all $|n| \geqslant 2$, we have
	\begin{equation*}
		\left\| \tilde{u}' \right\|_{L^1(n,n+1)} \leqslant 2^{\alpha} |n|^{-\alpha} \; ,
	\end{equation*}
	and consequently, $u_0$ satisfies Condition ($\mathcal{C}_{\mu, \alpha,2}$).
	\end{enumerate}
\end{5-REM4}

As above, we shall use several times the following basic lemma in the present section.

\begin{6-LEM2} \label{6-LEM2}
	Let $p \in [n,n+1]$, where $n \geqslant 1$ or $n \leqslant -2$. Then we have
	\begin{equation*}
		\frac{1}{2} \, |n| \leqslant |p| \leqslant 2 |n| \; .
	\end{equation*}
\end{6-LEM2}
\begin{proof}
	Firstly let us suppose that $n \geqslant 1$. Then
	\begin{equation*}
		\frac{1}{2} \, n \leqslant n \leqslant p \leqslant n + 1 \leqslant 2 \, n \; .
	\end{equation*}
	Now by supposing $n \leqslant -2$, we have
	\begin{equation*}
		\frac{1}{2} \, |n| = - \frac{1}{2} \, n \leqslant -(n+1) \leqslant -p = |p| \leqslant -n = |n| \leqslant 2 \, |n| \; .
	\end{equation*}
\end{proof}

To prove Theorem \ref{5-SCHRO3}, we start by splitting the infinite integration interval of the integral defining the solution formula \eqref{formula2} as follows: the singular frequency $0$ is the center of a sufficiently large but bounded interval, and we decompose the two remaining infinite intervals in an infinite union of disjoints bounded intervals. Thanks to that, the solution of the above evolution equation for an initial datum satisfying Condition (C3$_{\mu,\alpha,r}$) is actually a (infinite) sum of solutions of the same evolution equation but for initial data in bounded frequency bands. Then we follow the lines of the proof of Theorem \ref{4-SCHRO1} to apply the abstract results of Section \ref{Sec-Theory}, leading to a uniform estimate of each term of the sum in the cone $\mathfrak{C}_f(-\infty,+\infty)$, defined by the limits of $f'$ at infinity, as well as in its complement. Hence the series given by these uniform estimates provides a bound for the solution which is studied here. To assure the convergence of this series, we suppose that the decay at infinity of the Fourier transform of the initial datum is sufficiently fast as compared with the decay of the second derivative of the symbol.

\begin{5-SCHRO3} \label{5-SCHRO3}
	Suppose that the symbol $f$ satisfies Condition \emph{(S$_{\beta_+, \beta_-,R}$)} and that the initial datum $u_0$ satisfies Condition \emph{(C3$_{\mu, \alpha,r}$)}, where $\mu \in (0,1]$, $\alpha - \mu > \beta_-$ and $r \leqslant R$. Then we have
	\begin{equation*}
		\forall \, (t,x) \in \mathfrak{C}_f(-\infty,+\infty) \qquad \big| u_f(t,x) \big| \leqslant c^{(1)}(u_0,f) \, t^{- \frac{\mu}{2}} + c^{(2)}(u_0,f) \, t^{- \frac{1}{2}} \; ,
	\end{equation*}
	where the constants $c^{(1)}(u_0,f)$ and $c^{(2)}(u_0,f)$ are given by \eqref{C_1} and \eqref{C_(2)}, respectively. Moreover we have
	\begin{equation*}
		\forall \, (t,x) \in \mathfrak{C}_f(-\infty,+\infty)^c \qquad \big| u_f(t,x) \big| \leqslant c_c^{(1)}(u_0,f) \, t^{- \mu} + c_c^{(2)}(u_0,f) \, t^{-1} \; ,
	\end{equation*}
	where the constants $c_c^{(1)}(u_0,f)$ and $c_c^{(2)}(u_0,f)$ are given by \eqref{C_(1)} and \eqref{C_2}, respectively.\\
	The space-time cone $\mathfrak{C}_f(-\infty,+\infty)$ is defined by
	\begin{equation*}
		\mathfrak{C}_f(-\infty,+\infty) := \left\{ (t,x) \in (0,+\infty) \times \R \, \Big| \, a < \frac{x}{t} < b \right\} \; ,
	\end{equation*}
	where $\displaystyle a := \lim_{p \rightarrow - \infty} f'(p)$ and $\displaystyle b := \lim_{p \rightarrow + \infty} f'(p)$.
\end{5-SCHRO3}

\begin{proof}
	We recall that the solution of the initial value problem \eqref{evoleq} can be written as follows,
	\begin{equation*}
		\forall \, (t,x) \in (0,+\infty) \times \R \qquad u_f(t,x) = \int_{\R} U(p) \, e^{i t \psi(p)} \, dp \; ,
	\end{equation*}
	where
	\begin{equation*}
		\left\{ \begin{array}{rl}
				& \displaystyle \forall \, p \in \R \backslash \{0 \} \qquad U(p) := \frac{1}{2\pi}  \, \tf u_0(p) = \frac{1}{2\pi}  \, |p|^{\mu-1} \, \tilde{u}(p) \; , \\
				& \vspace{-0.3cm} \\
				& \displaystyle \forall \, p \in \R \qquad \psi(p) := \frac{x}{t} \, p - f(p) \; .
		\end{array} \right.
	\end{equation*}
	Let us define $N \in \mathbb{N}$ and $\mathfrak{S}_N \subseteq \mathbb{Z}$ as follows,
		\begin{equation*}
			N := \lceil R \rceil + 1 \qquad , \qquad \mathfrak{S}_N = \mathbb{Z} \backslash \{-N,\dots, N-1 \} \; ,
		\end{equation*}
		where $\lceil . \rceil$ is the ceiling function, and let us split the integral,
		\begin{align*}
			\int_{\R} U(p) e^{it \psi(p)} \, dp	& = \int_{\R} \chi_{[-N,N)}(p) U(p) \, e^{it \psi(p)} \, dp + \int_{\R} \sum_{n \in \mathfrak{S}_N} \chi_{[n,n+1)}(p) U(p) \, e^{it \psi(p)} \, dp \\
												& = \int_{-N}^{N} U(p) \, e^{it \psi(p)} \, dp + \sum_{n \in \mathfrak{S}_N} \int_n^{n+1} U(p) \, e^{it \psi(p)} \, dp \; ,
		\end{align*}
		where $\chi_{[n,n+1)}$ is the characteristic function of the interval $[n,n+1)$. Now we divide the proof into two parts with respect to the value of $\frac{x}{t}$.
		\begin{enumerate}
		\item \textit{Case $\frac{x}{t} \in (a,b)$.} In this case, $\frac{x}{t}$ belongs to $(a,b)$, that is to say it belongs to $f'(\R)$. Therefore the phase $\psi$ has a unique stationary point which belongs to $\R$.\\
		To estimate the integral on $[-N,N]$ in this case, we apply Theorem \ref{3-COR1} for $\rho = 2$ on $[-N,0]$ and on $[0,N]$ by following the lines of the proof of Theorem \ref{4-SCHRO1} in the case i) which gives
		\begin{align}
			\forall \, (t,x) \in \mathfrak{C}_f(-\infty,+\infty) \qquad \left| \int_{-N}^N U(p) \, e^{i t \psi(p)} \, dp \right|	& \leqslant \left| \int_{-N}^{0} \dots \right| + \left| \int_0^N \dots \right| \nonumber \\
						& \leqslant \Big(c_{-N}^{(1)}(u_0,f) + c_{+N}^{(1)}(u_0,f) \Big) \, t^{-\frac{\mu}{2}} \nonumber \\
						& \label{C_1} =: c^{(1)}(u_0,f) \, t^{-\frac{\mu}{2}} \; ,
		\end{align}
		with
		\begin{align*}
			& \bullet \quad c_{-N}^{(1)} (u_0,f) := \frac{1}{2 \pi} \, \frac{3}{\mu} \left\| \tilde{u} \right\|_{L^{\infty}(-N,0)} + \frac{1}{\pi} \, \Big( 4 \left\| \tilde{u} \right\|_{L^{\infty}(-N,0)} + \left\| \tilde{u}' \right\|_{L^1(-N,0)} \Big) \, m_{-N}^{\; -1} \; , \\
			& \bullet \quad c_{+N}^{(1)} (u_0,f) := \frac{1}{2 \pi} \, \frac{3}{\mu} \left\| \tilde{u} \right\|_{L^{\infty}(0,N)} + \frac{1}{\pi} \, \Big( 4 \left\| \tilde{u} \right\|_{L^{\infty}(0,N)} + \left\| \tilde{u}' \right\|_{L^1(0,N)} \Big) \, m_{+N}^{\; -1} \; ,
		\end{align*}
		and
		\begin{equation*}
			m_{- N} := \min_{p \in [-N, 0]} f''(p) > 0 \qquad , \qquad m_{+ N} := \min_{p \in [0,N]} f''(p) > 0 \; .
		\end{equation*}
		Now let us study each term of the series. By hypothesis, $U$ has no singular points in $[n,n+1]$ for $n \in \mathfrak{S}_N$. As above, we can apply Theorem \ref{3-COR1} for $\rho = 2$ and $\mu = 1$, and we obtain
		\begin{equation*}
			\forall \, (t,x) \in \mathfrak{C}_f(-\infty,+\infty) \qquad \left| \int_{n}^{n+1} U(p) \, e^{it \psi(p)} \, dp \right| \leqslant c_n^{(2)}(u_0,f) \, t^{-\frac{1}{2}} \; ;
		\end{equation*}
		the constant $c_n^{(2)}(u_0,f) > 0$ is given by
		\begin{equation*}
			c_n^{(2)}(u_0,f)	:= \frac{1}{\pi} \, \left\| U \right\|_{L^{\infty}(n,n+1)} + \frac{1}{\pi} \, \Big( 3 \left\| U \right\|_{L^{\infty}(n,n+1)} + \left\| U' \right\|_{L^1(n,n+1)} \Big) \, m_{n}^{\; -1} \; ,
		\end{equation*}
		with $\displaystyle m_{n} := \min_{p \in [n,n+1]} f''(p) > 0$.\\
		The following step is to prove the summability of the sequence $\big\{ c_n^{(2)}(u_0,f) \big\}_{n \in \mathfrak{S}_N}$. On the one hand, we have by using the hypothesis on $u_0$ and Lemma \ref{6-LEM2},
		\begin{equation*}
			\left\| U \right\|_{L^{\infty}(n,n+1)} \leqslant 2^{1-\mu} |n|^{\mu-1} \, M 2^{\alpha} |n|^{-\alpha} = 2^{1-\mu + \alpha} M \, |n|^{\mu-1-\alpha} \; .
		\end{equation*}
		Moreover
		\begin{align}
			\left\| U' \right\|_{L^1(n,n+1)}	& \leqslant \int_n^{n+1} (1-\mu) |p|^{\mu-2} \left| \tilde{u} (p) \right| dp + \int_n^{n+1} |p|^{\mu-1} \left| \tilde{u}' (p) \right| dp \nonumber \\
												& \leqslant \left\| \tilde{u} \right\|_{L^{\infty}(n,n+1)} \int_n^{n+1} (1-\mu) |p|^{\mu-2} dp + 2^{1-\mu} \, |n|^{\mu-1} \left\| \tilde{u}' \right\|_{L^1(n,n+1)} \nonumber \\
												& \label{estderivu} \leqslant M 2^{\alpha} |n|^{-\alpha} \, 2^{1-\mu} |n|^{\mu-1} + 2^{1-\mu} |n|^{\mu-1} \, M' |n|^{-\alpha} \\
												& = 2^{1-\mu} \left( 2^{\alpha} M + M' \right) |n|^{\mu-1-\alpha} \nonumber \; ,
		\end{align}
		where the hypothesis $\big\| \tilde{u}' \big\|_{L^1(n,n+1)} \leqslant M' |n|^{-\alpha}$ was used to get \eqref{estderivu}. On the other hand, we have by the hypothesis on the symbol $f$,
		\begin{equation*}
			f''(p) \geqslant c_- \, |p|^{-\beta_-} \geqslant c_- \, 2^{-\beta_-} |n|^{-\beta_-} \; .
		\end{equation*}
		It follows
		\begin{equation*}
			m_n^{\; -1} \leqslant \frac{2^{\beta_-}}{c_-} \, |n|^{\beta_-} \; .
		\end{equation*}
		Then we obtain
		\begin{align*}
			c_n^{(2)}(u_0,f)	& = \frac{1}{\pi} \, \left\| U \right\|_{L^{\infty}(n,n+1)} + \frac{1}{\pi} \, \Big( 3 \left\| U \right\|_{L^{\infty}(n,n+1)} + \left\| U' \right\|_{L^1(n,n+1)} \Big) \, m_{n}^{\; -1} \\
						& \leqslant \frac{2^{1-\mu +\alpha} M}{\pi} \, |n|^{\mu-1-\alpha} + 3 \, \frac{2^{1-\mu+\alpha+\beta_-} M}{\pi \, c_-} \, |n|^{\mu-1-\alpha+\beta_-} \\
						& \qquad \qquad + \, \frac{2^{1-\mu+\beta_-}  \left( 2^{\alpha} M + M' \right)}{\pi \, c_-} \, |n|^{\mu-1-\alpha+\beta_-} \; .
		\end{align*}
		Since $\alpha - \mu > \beta_-$, the sequence $\big\{ c_n^{(2)}(u_0,f) \big\}_{n \in \mathfrak{S}_N}$ is summable. It follows
		\begin{equation*}
			\left| \sum_{n \in \mathfrak{S}_N} \int_n^{n+1} U(p) \, e^{it \psi(p)} \, dp \right|	\leqslant \sum_{n \in \mathfrak{S}_N} \left| \int_n^{n+1} U(p) \, e^{it \psi(p)} \, dp \right| \leqslant \left( \sum_{n \in \mathfrak{S}_N} c_n^{(2)}(u_0,f) \right) t^{-\frac{1}{2}} \; .
		\end{equation*}
		Then it is possible to bound the last series by employing the following estimate of the Riemann Zeta function,
		\begin{equation*}
			\forall \, \sigma > 1 \qquad \sum_{n \in \N^*} n^{-\sigma} \leqslant \frac{\sigma}{\sigma - 1} \; .
		\end{equation*}
		Hence
		\begin{align}
			\sum_{n \in \mathfrak{S}_N} c_n^{(2)}(u_0,f)	& \leqslant \frac{2^{2 - \mu + \alpha} M}{\pi} \, \frac{\alpha + 1 - \mu}{\alpha - \mu} \, + 3 \, \frac{2^{2-\mu+\alpha+\beta_-} M}{\pi \, c_-} \, \frac{\alpha + 1 - \mu - \beta_-}{\alpha - \mu - \beta_-} \nonumber \\
													& \qquad + \, \frac{2^{2-\mu+\beta_-}  \left( 2^{\alpha} M + M' \right)}{\pi \, c_-} \, \frac{\alpha + 1 - \mu - \beta_-}{\alpha - \mu - \beta_-} \nonumber \\
													& = \frac{2^{2 - \mu + \alpha} M}{\pi} \, \frac{\alpha + 1 - \mu}{\alpha - \mu} \, + \, \frac{2^{2-\mu+\beta_-} ( 2^{\alpha+2} M + M' )}{\pi \, c_-} \, \frac{\alpha + 1 - \mu - \beta_-}{\alpha - \mu - \beta_-} \nonumber \\
													& =: c^{(2)}(u_0,f) \; . \label{C_(2)}
		\end{align}
		Hence we obtain finally for all $(t,x) \in \mathfrak{C}_f(-\infty,+\infty)$,
		\begin{equation*}
			\big| u_f(t,x) \big| \leqslant c^{(1)}(u_0,f) \, t^{- \frac{\mu}{2}} + c^{(2)}(u_0,f) \, t^{- \frac{1}{2}} \; .
		\end{equation*}
		\item \textit{Case $\frac{x}{t} \notin (a,b)$.}	We note that $\psi$ has no stationary points on $\R$ and that $(t,x)$ belongs to $\mathfrak{C}_f(-\infty,+\infty)^c$ in this case.\\
		We start by estimating the integral on $[-N,N]$. To do so, we apply Theorem \ref{3-VDC4} on $[-N,0]$ and on $[0,N]$, providing
	\begin{align}
		\forall \, (t,x) \in \mathfrak{C}_f(-\infty,+\infty)^c \qquad \left| \int_{-N}^N U(p) \, e^{i t \psi(p)} \, dp \right|	& \leqslant \left| \int_{-N}^{0} \dots \right| + \left| \int_0^N \dots \right| \nonumber \\
					& \leqslant \Big(\tilde{c}_{-N}^{(1)}(u_0,f) + \tilde{c}_{+N}^{(1)}(u_0,f) \Big) \, t^{-\mu} \nonumber \\
						& =: c_c^{(1)}(u_0,f) \, t^{-\mu} \; , \label{C_(1)}
	\end{align}
	with
	\begin{align*}
		& \bullet \quad \tilde{c}_{-N}^{(1)}(u_0,f) := \frac{1}{2 \pi} \, \frac{1}{\mu} \left\| \tilde{u} \right\|_{L^{\infty}(-N,0)} + \frac{1}{2\pi} \, \Big( 4 \left\| \tilde{u} \right\|_{L^{\infty}(-N,0)} + \left\| \tilde{u}' \right\|_{L^1(-N,0)} \Big) \, \tilde{m}_{-N}^{\; -1} \; , \\
		& \bullet \quad \tilde{c}_{+N}^{(1)}(u_0,f) := \frac{1}{2 \pi} \, \frac{1}{\mu} \left\| \tilde{u} \right\|_{L^{\infty}(0,N)} + \frac{1}{2\pi} \, \Big( 4 \left\| \tilde{u} \right\|_{L^{\infty}(0,N)} + \left\| \tilde{u}' \right\|_{L^1(0,N)} \Big) \, \tilde{m}_{+N}^{\; -1} \; .
	\end{align*}
	The terms $\tilde{m}_{+N}, \tilde{m}_{-N} > 0$ are defined as follows,
	\begin{align*}
		& \bullet \quad \forall \, p \in [-N,0] \quad \left| \psi'(p) \right| = \left| \frac{x}{t} - f'(p) \right| \geqslant \min \big\{ f'(-N) - a , b - f'(0) \big\} =: \tilde{m}_{-N} \; , \\
		& \bullet \quad \forall \, p \in [0,N] \quad \left| \psi'(p) \right| = \left| \frac{x}{t} - f'(p) \right| \geqslant \min \big\{ f'(0) - a , b - f'(N) \big\} =: \tilde{m}_{+N} \; .
	\end{align*}
	Now we study the terms of the series. The amplitude $U$ has no singular points in $[n,n+1]$ for $n \in \mathfrak{S}_N$ by hypothesis, so Theorem \ref{3-VDC4} is applicable once again on the interval $[n,n+1]$ with $\mu = 1$ and it gives
	\begin{equation*}
		\forall \, (t,x) \in \mathfrak{C}_f(-\infty,+\infty)^c \qquad \left| \int_{n}^{n+1} U(p) \, e^{it \psi(p)} \, dp \right| \leqslant \tilde{c}_n^{(2)}(u_0,f) \, t^{-1} \; ;
	\end{equation*}
	the constant $\tilde{c}_n^{(2)}(u_0,f) > 0$ is defined by
	\begin{equation*}
		\tilde{c}_n^{(2)}(u_0,f) := \frac{1}{2 \pi} \Big( 3 \left\| U \right\|_{L^{\infty}(n,n+1)} + \left\| U' \right\|_{L^1(n,n+1)} \Big) \, \tilde{m}_n^{\; -1} \; ,
	\end{equation*}
	with $\displaystyle \tilde{m}_n := \min \big\{ f'(n) - a , b - f'(n+1) \big\} > 0$. As in the previous case, by using the fact that $u_0$ satisfies Condition (C3$_{\mu,\alpha,r}$) and using Lemma \ref{6-LEM2}, one can show that
	\begin{equation*}
		\left\| U \right\|_{L^{\infty}(n,n+1)} \leqslant 2^{1-\mu + \alpha} M \, |n|^{\mu-1-\alpha} \; ,
	\end{equation*}
	and
	\begin{equation*}
		\left\| U' \right\|_{L^1(n,n+1)} \leqslant 2^{1-\mu} \left( 2^{\alpha} M + M' \right) |n|^{\mu-1-\alpha} \nonumber \; .
	\end{equation*}
	Furthermore the point ii) of Lemma \ref{5-LEM3} implies
	\begin{equation*}
		\tilde{m}_n \geqslant \frac{c_-}{\beta_- - 1} \, \min \left\{ |n|^{1-\beta_-} , |n+1|^{1-\beta_-} \right\} \geqslant \frac{c_-}{\beta_- - 1} \, 2^{1-\beta_-} |n|^{1-\beta_-} \; ,
	\end{equation*}
	where we used Lemma \ref{6-LEM2} one more time. Then we obtain
	\begin{align*}
		\tilde{c}_n^{(2)}(u_0,f)	& = \frac{1}{2 \pi} \Big( 3 \left\| U \right\|_{L^{\infty}(n,n+1)} + \left\| U' \right\|_{L^1(n,n+1)} \Big) \, \tilde{m}_n^{\; -1} \\
					& \leqslant \frac{\beta_- - 1}{2 \pi c_-} \Big( 3 \times 2^{-\mu+\alpha+\beta_-} \, M \, + 2^{-\mu + \beta_-} \left(2^{\alpha} M + M' \right) \Big) \, |n|^{\mu-2-\alpha+\beta_-} \; .
	\end{align*}
	Then the summability of the sequence $\big\{ \tilde{c}_n^{(2)}(u_0,f) \big\}_{n \in \mathfrak{S}_N}$ follows from the assumption $\alpha - \mu > \beta_- > \beta_- - 1$, and we have
	\begin{align}
		\sum_{n \in \mathfrak{S}_N} \tilde{c}_n^{(2)}(u_0,f) 	& \leqslant \frac{\beta_- - 1}{ \pi c_-} \Big( 3 \times 2^{-\mu+\alpha+\beta_-} \, M \, + 2^{-\mu + \beta_-} (2^{\alpha} M + M') \Big) \, \frac{\alpha + 2 - \mu - \beta_-}{\alpha + 1 - \mu - \beta_-} \nonumber \\
													& =: c_c^{(2)}(u_0,f) \; . \label{C_2}
	\end{align}
	It follows that
	\begin{equation*}
		\left| \sum_{n \in \mathfrak{S}_N} \int_n^{n+1} U(p) \, e^{it \psi(p)} \, dp \right| \leqslant \left( \sum_{n \in \mathfrak{S}_N} \tilde{c}_n^{(2)}(u_0,f) \right) t^{-1} \leqslant c_c^{(2)}(u_0,f) \, t^{-1} \; .
	\end{equation*}
	We obtain finally for all $(t,x) \in \mathfrak{C}_f(-\infty,+\infty)^c$,
	\begin{equation*}
		\big| u_f(t,x) \big| \leqslant c_c^{(1)}(u_0,f) \, t^{-\mu} + c_c^{(2)}(u_0,f) \, t^{-1} \; .
	\end{equation*}
	\end{enumerate}
\end{proof}

\vspace{4mm}

Let us illustrate this intrinsic localization phenomenon by applying Theorem \ref{5-SCHRO3} to the solution formula of the Klein-Gordon equation. In particular, this example shows the natural appearance of equations of type \eqref{evoleq} when studying higher order in time hyperbolic equations.\\
The Klein-Gordon equation on $\R$ is defined as follows:
\begin{equation} \label{eqkg}
	\left\{ \begin{array}{rl}
			& \hspace{-2mm} \big[ \, \partial_{tt} - c^2 \partial_{xx} + c^4 \, \big] u_{KG}(t) = 0 \\ [2mm]
			& \hspace{-2mm} u_{KG}(0) = u_0 \quad , \quad \partial_t \, u_{KG}(0) = v_0
	\end{array} \right. \; ,
\end{equation}
for $t \geqslant 0$, where $c > 0$ is a constant. In terms of quantum mechanics, the constant $c$ represents the speed of light and the solution is the wave function of a spinless relativistic free particle with mass $m = 1$. By assuming that $u_0, v_0 \in \mathcal{S}'(\R)$, one can furnish a solution formula which belongs to $\displaystyle \mathcal{C}^{2}\big( \R_+ , \mathcal{S}'(\R) \big)$,
\begin{equation} \label{KGformula}
	u_{KG}(t) = \tf^{-1} \Big( e^{-i t f_{KG}} a_+(u_0,v_0) \Big) + \tf^{-1} \Big( e^{i t f_{KG}} a_-(u_0,v_0) \Big) =: u_+(t) + u_-(t) \; ,
\end{equation}
where the symbol $f_{KG}: \R \longrightarrow \R$ is given by $f_{KG}(p) = \sqrt{c^4 + c^2 \, p^2 \,}$, and the tempered distributions $a_+(u_0,v_0)$ and $a_-(u_0,v_0)$ are defined by
\begin{equation*} 
	a_+(u_0,v_0) := \frac{1}{2} \left( \tf u_0 + \frac{i}{f_{KG}} \, \tf v_0 \right) \quad , \quad a_-(u_0,v_0) := \frac{1}{2} \left( \tf u_0 - \frac{i}{f_{KG}} \, \tf v_0 \right) \; .
\end{equation*}
In this context, we note that $u_+$ and $u_-$ solve respectively the equations
\begin{equation} \label{eqpm}
	\left\{ \begin{array}{rl}
			& \hspace{-2mm} \left[ i \, \partial_t - f_{KG} \big( D \big) \right] u_+(t) = 0 \\ [2mm]
			& \hspace{-2mm} u_+(0) = \tf^{-1} a_+(u_0,v_0)
	\end{array} \right. \quad , \qquad \left\{ \begin{array}{rl}
			& \hspace{-2mm} \left[ i \, \partial_t + f_{KG} \big( D \big) \right] u_-(t) = 0 \\ [2mm]
			& \hspace{-2mm} u_-(0) = \tf^{-1} a_-(u_0,v_0)
	\end{array} \right. \; ,
\end{equation}
for $t \geqslant 0$. In particular, the solution formula \eqref{KGformula} defines a complex-valued function on $(0,+\infty) \times \R$ when $\tf^{-1} a_+(u_0,v_0)$ and $\tf^{-1} a_-(u_0,v_0)$ satisfy Condition (C3$_{\mu, \alpha,r}$). In the following result, we provide estimates of the solution of the Klein-Gordon equation \eqref{eqkg} coming from estimates for the evolution equations given in \eqref{eqpm}. The proof consists mainly in showing that the symbol $f_{KG}$ satisfies Condition (S$_{\beta_+,\beta_-,R}$), for certain $\beta_+, \beta_-, R$. Theorem \ref{5-SCHRO3} is then applicable, and the resulting estimates indicate that the solution of the Klein-Gordon equation \eqref{eqkg} is time-asymptotically localized in the space-time cone
\begin{equation}
	\mathfrak{C}_{f_{KG}}(-\infty,+\infty) = \left\{ (t,x) \in (0,+\infty) \times \R \, \Big| \, -c < \frac{x}{t} < c \right\} \; ,
\end{equation}
which is actually the light cone issued by the origin.

\begin{5-COR3} \label{5-COR3}
	Let $u_{KG}: \R_+ \times \R \longrightarrow \C$ be the solution of the Klein-Gordon equation on $\R$ with $u_0, v_0 \in \mathcal{S}'(\R)$ such that $\tf^{-1} a_+(u_0,v_0)$ and $\tf^{-1} a_-(u_0,v_0)$ satisfy Condition \emph{(C3$_{\mu, \alpha,r}$)}, with $\mu \in (0,1]$, $\alpha - \mu > 3$ and $r \leqslant c$. Then we have
	\begin{equation*}
		\forall \, (t,x) \in \mathfrak{C}_{f_{KG}}(-\infty,+\infty) \qquad \big| u_{KG}(t,x) \big| \leqslant c^{(1)}(u_0, v_0, f_{KG}) \, t^{- \frac{\mu}{2}} + c^{(2)}(u_0, v_0, f_{KG}) \, t^{- \frac{1}{2}} \; ,
	\end{equation*}
	and
	\begin{equation*}
		\forall \, (t,x) \in \mathfrak{C}_{f_{KG}}(-\infty,+\infty)^c \qquad \big| u_{KG}(t,x) \big| \leqslant c_c^{(1)}(u_0, v_0, f_{KG}) \, t^{- \mu} + c_c^{(2)}(u_0, v_0, f_{KG}) \, t^{-1} \; .
	\end{equation*}
	All the constants can be computed from Theorem \ref{5-SCHRO3}.
\end{5-COR3}

\begin{proof}
	First of all, let us remark that one can follow the lines of the proof of Theorem \ref{5-SCHRO3} to establish very similar estimates for the solution of the evolution equation \eqref{evoleq} when $-f$ satisfies Condition (S$_{\beta_+,\beta_-,R}$), that is to say when $f'' < 0$. In this case, $a$ is the limit of $f'$ at $+ \infty$ and $b$ the limit of $f'$ at $-\infty$. In the present proof, this remark assures that Theorem \ref{5-SCHRO3} is applicable to both equations given in \eqref{eqpm} if the symbol $f_{KG}$ verifies Condition (S$_{\beta_+,\beta_-,R}$).\\
	Now we provide the first and the second derivative of $f_{KG}$,
	\begin{equation*}
		\forall \, p \in \R \qquad (f_{KG})'(p) = \frac{c \, p}{\sqrt{c^2 + p^2}} \qquad , \qquad (f_{KG})''(p) = c^3 \left( \frac{c^2}{p^2} \, + 1 \right)^{-\frac{3}{2}} \, |p|^{-3} \; .
	\end{equation*}
	By noting that the following inequalities are true,
	\begin{equation*}
		\forall \, |p| \geqslant c \qquad 2^{-\frac{3}{2}} \, c^3 \leqslant c^3 \left( \frac{c^2}{p^2} \, + 1 \right)^{-\frac{3}{2}} \leqslant c^3 \; ,
	\end{equation*}
	we deduce that $f_{KG}$ satisfies Condition (S$_{3,3,c}$). Moreover one can see that the limits of $(f_{KG})'$ at $-\infty$ and $+\infty$ are given by $-c$ and $c$ respectively. It follows that Theorem \ref{5-SCHRO3} is applicable to the solutions of the equations \eqref{eqpm}, furnishing the estimates of the solution $u_{KG}$ inside the cone $\mathfrak{C}_{f_{KG}}(-\infty,+\infty)$ and outside.
\end{proof}

\begin{5-REM1} \label{reduction} \em
	As explained above, the Klein-Gordon equation furnishes a simple setting illustrating the localization phenomenon exhibited in Theorem \ref{5-SCHRO3}. Nevertheless it is also possible to consider more general equations, as for example:
\begin{equation} \label{eqhyp}
	\left\{ \begin{array}{rl}
			& \displaystyle \hspace{-2mm} \left[ \partial_t^n - \sum_{k = 0}^n (-i)^{k-n} \, a_k \, \partial_x^k \right] \hspace{-1mm} u(t) = 0 \\ [5mm]
			& \hspace{-2mm} u(0) = u_0 \; , \; \partial_t u(0) = u_1 \; , \; ... \; , \partial_t^{n-1} u(0) = u_{n-1}
	\end{array} \right. \; ,
\end{equation}
where $t \geqslant 0$, $n \geqslant 2$ is an even number and $\{a_k\}_{k=0,...,n}$ is a set of real numbers such that
\begin{equation*}
	\forall \, p \in \R \qquad F(p) := \sum_{k = 0}^n a_k \, p^k > 0 \; .
\end{equation*} 
Applying the Fourier transform to \eqref{eqhyp}, resolving the resulting ODE and applying the inverse Fourier transform show that the solution formula of equation \eqref{eqhyp} is a sum of $n$ terms such that two of them are solutions of equations of type \eqref{evoleq} with symbols $f_n$ and $-f_n$, where $f_n := \sqrt[n]{F \,}$, respectively. Under suitable hypotheses on the coefficients $a_k$, the symbol $f_n$ satisfies Condition (S$_{\beta_+, \beta_-,R}$), for certain $\beta_+, \beta_-$ and $R$, implying the above time-asymptotic localization for the two mentioned terms.
\end{5-REM1}

\vspace{0.5cm}

\noindent \textbf{Acknowledgements:}\\
The author thanks E. Creusé for valuable support and R. Haller-Dintelmann for valuable discussions. The author thanks also F. Ali Mehmeti for the helpful and numerous discussions which permit to improve the content of this paper.\\
The author has been supported by a research grant from the excellence laboratory in mathematics and physics CEMPI and the region Nord-Pas-de-Calais (France).

\end{document}